\documentclass[final]{siamltex}
\usepackage{amsmath,amsfonts,amssymb}
\usepackage{etoolbox}
\usepackage{latexsym}
\usepackage{epsfig,overpic}
\usepackage{comment}
\usepackage{todonotes}
\usepackage{pgfplots}
\usepackage{tikz}
\usepackage{microtype}
\usepackage[colorlinks=true,citecolor=red,linkcolor=blue,pdfborder={0 0 0}]{hyperref}
\usetikzlibrary{shapes,arrows,snakes,calendar,matrix,backgrounds,folding,calc,positioning}
\DeclareGraphicsExtensions{.pdf,.eps,.ps,.eps.gz,.ps.gz,.eps.Y}
\usepackage{color,colortbl}
\usepackage{tikz}

\definecolor{mildblue}{HTML}{0000b0}
\definecolor{mildpurple}{HTML}{400090}
\definecolor{mildgreen}{HTML}{007000}
\definecolor{mildred}{HTML}{b00000}
\definecolor{lightblue}{HTML}{57c5d9}
\definecolor{lightred}{HTML}{e35555}
\definecolor{grey}{HTML}{BBBBBB}
\definecolor{lightgrey}{HTML}{EEEEEE}
\definecolor{black}{HTML}{000000}
\definecolor{myblockcolor}{HTML}{002e6e}

\newcommand{\bb}{\mathbf b}
\newcommand{\bc}{\mathbf{c}}
\newcommand{\bd}{\mathbf{d}}
\newcommand{\be}{\mathbf e}

\newcommand{\bn}{\mathbf n}
\newcommand{\bp}{\mathbf p}
\newcommand{\bq}{\mathbf{q}}

\newcommand{\bs}{\mathbf s}
\newcommand{\bu}{\mathbf u}
\newcommand{\bv}{\mathbf v}
\newcommand{\bw}{\mathbf w}
\newcommand{\bx}{\mathbf x}
\newcommand{\by}{\mathbf y}
\newcommand{\bz}{\mathbf z}
\newcommand{\bX}{\mathbf X}

\newcommand{\ccg}{\cellcolor{lightgrey}}

\newcommand{\bng}{{\mathbf n}_{\mathcal{G}}}

\newcommand{\pbx}[1]{{\color{mildgreen} {\mathbf x}^{#1}}}
\newcommand{\pby}[1]{{\color{mildpurple} {\mathbf y}^{#1}}}
\newcommand{\pbb}[1]{{\color{mildpurple} {\mathbf b}^{#1}}}

\newcommand{\pbxx}[1]{{\color{mildpurple} {\mathbf x}^{#1}}}
\newcommand{\pbc}[1]{{\color{mildgreen} {\mathbf c}^{#1}}}
\newcommand{\pbd}[1]{{\color{mildpurple} {\mathbf d}^{#1}}}

\newcommand{\pbbu}[1]{{\color{mildred} { \underline{\mathbf u}}^{#1}}}
\newcommand{\pbbv}[1]{{\color{mildgreen} { \underline{\mathbf v}}^{#1}}}
\newcommand{\pbbw}[1]{{\color{mildblue} {\underline{\mathbf w}}^{#1}}}

\newcommand{\pbbx}[1]{{\color{mildred} {\underline{\mathbf x}}^{#1}}}
\newcommand{\pbbc}[1]{{\color{mildgreen} {\underline{\mathbf c}}^{#1}}}
\newcommand{\pbbd}[1]{{\color{mildblue} {\underline{\mathbf d}}^{#1}}}

\newcommand{\pbbe}[2]{{\color{yellow!40!black} {\underline{\mathbf #1}}^{#2}}}
\newcommand{\pbbf}[2]{{\color{orange!60!black} {\underline{\mathbf #1}}^{#2}}}
\newcommand{\pbbg}[2]{{\color{purple!60!black} {\underline{\mathbf #1}}^{#2}}}

\newcommand{\bconv}{\mathbf{conv}}

\newcommand{\mI}{\mathcal{I}}
\newcommand{\mH}{\mathcal{H}}
\newcommand{\mD}{\mathcal{D}}

\newcommand{\mP}{\mathcal{P}}
\newcommand{\mQ}{\mathcal{Q}}
\newcommand{\mB}{\mathcal{B}}
\newcommand{\mT}{\mathcal{T}}
\newcommand{\mG}{\mathcal{G}}

\newcommand{\T}{\mathcal T}

\newcommand{\Div}{\mathop{\rm div}}

\newcommand{\rr}{\mathbb{R}}

\newcommand{\zz}{\mathbb{Z}}
\newcommand{\nn}{\mathbb{N}}

\newcommand{\Gs}{{\Gamma_\ast}}
\newcommand{\Gsh}{{\Gamma_{\ast,h}}}

\renewcommand{\div}{\textrm{div}\ \!}

\newcommand{\ms}{m_s}
\newcommand{\mt}{m_t}
\newcommand{\ns}{n_s}
\newcommand{\nt}{n_t}
\newcommand{\eocs}{eoc$_s$}
\newcommand{\eoct}{eoc$_t$}
\newcommand{\pr}{\%}

\newtheorem{remark}{Remark}[section]

\newtoggle{extended}
\toggletrue{extended}

\begin{document}

\title{The Nitsche XFEM-DG space-time method and its implementation in three space dimensions}
\author{Christoph Lehrenfeld
\thanks{Institut f\"ur Geometrie und Praktische Mathematik, RWTH Aachen,
D-52056 Aachen, Germany; email: {\tt lehrenfeld@igpm.rwth-aachen.de}}}
\maketitle
\begin{abstract}
In the recent paper [C. Lehrenfeld, A. Reusken, \emph{SIAM J. Num. Anal.}, 51 (2013)] a new finite element discretization method for a class of two-phase mass transport problems is presented and analyzed. The transport problem describes mass transport in a domain with an evolving interface. Across the evolving interface a jump condition has to be satisfies. The discretization in that paper is a space-time approach which combines a discontinuous Galerkin (DG) technique (in time) with an extended finite element method (XFEM). Using the Nitsche method the jump condition is enforced in a weak sense. 
While the emphasis in that paper was on the analysis and one dimensional numerical experiments the main contribution of this paper is the discussion of implementation aspects for the spatially three dimensional case. 
As the space-time interface is typically given only implicitly as the zero-level of a level-set function, we construct a piecewise planar approximation of the space-time interface. This discrete interface is used to divide the space-time domain into its subdomains. An important component within this decomposition is a new method for dividing four-dimensional prisms intersected by a piecewise planar space-time interface into simplices. Such a subdivision algorithm is necessary for numerical integration on the subdomains as well as on the space-time interface. These numerical integrations are needed in the implementation of the Nitsche XFEM-DG method in three space dimensions. Corresponding numerical studies are presented and discussed.
\end{abstract}
\begin{keywords}finite elements, evolving surface, parabolic PDE, space-time, two-phase flow, computations, XFEM, 4D quadrature, discontinuous galerkin
\end{keywords}
\begin{AMS}
65D30, 65M30
\end{AMS}
\section[Introduction]{Introduction}\label{sec:introduction}
We consider two immiscible incompressible fluids which are contained in a convex polygonal domain $\Omega \subset \rr^d$, $d=2,3$. 
We assume (for simplicity) that one fluid is completely surrounded by the other.
A typical configuration of this sort is a droplet moving in a flow field. 
The domain containing the surrounding fluid is $\Omega_2$, whereas the inner phase is $\Omega_1 = \Omega \setminus \Omega_2$. 
Both domains are separated by a smooth interface $\Gamma := \bar{\Omega}_1 \cap \bar{\Omega}_2$. 
The hydrodynamic evolution of the fluids is typically modeled by the incompressible Navier-Stokes equations. 
Suitable interface conditions have to be added to account for continuity and momentum balance which includes surface tension forces. 
In \cite{Bothe04,Gross04,Muradoglu07,GReusken2011} one can find a more detailed discussion on the model.

For the mass transport problem under consideration in this paper we assume that the velocity field $\bw$ as well as the evolution of the domains $\Omega_1(t)$, $\Omega_1(t)$ (and the interface $\Gamma(t)$) are known and there holds $\div \bw = 0$ and $V_{\Gamma} = \bw \cdot \bn$ on $\Gamma(t)$ with $V_{\Gamma}$ the velocity of the interface motion in normal direction and $\bn$ the unit normal of $\partial \Omega_1(t)$. 
The latter assumption states that the relative convection velocity across the moving interface $\Gamma(t)$ is zero. 

We consider the transport of a dissolved species $u$ in such a two-phase flow configuration. A standard model consists of convection-diffusion equations within both domains, where Fick's law is applied as the diffusion model:
\begin{subequations}
  \begin{equation}
    \frac{\partial u}{\partial t} + \bw \cdot \nabla u - \Div (\alpha\nabla u)  = f \quad \text{in}~~ \Omega_m(t), ~m=1,2,~t \in [0,T], \label{eq1} \\
  \end{equation}
with $f$ a source term and $\alpha= \alpha(\bx,t)$ the diffusion coefficient which is assumed to be piecewise constant,
\[ 
 \alpha= \alpha_m(t) >0 \quad \text{in}~~\Omega_m(t).
\]
and different in the subdomains $\alpha_1 \neq \alpha_2$. The solutions in the subdomains are coupled through interface conditions
  \begin{eqnarray}
    {[}\alpha \nabla u \cdot \bn{]}_{\Gamma}  & =& 0 \quad \text{at}~~\Gamma(t) , \label{eq2} \\
    {[}\beta u{]}_{\Gamma} & = &0 \quad \text{at}~~\Gamma(t), \label{eq3}
  \end{eqnarray}
where $[v]_{\Gamma}$ denotes the jump across the interface $\Gamma$ of a (sufficiently smooth) quantity $v$, $[v]_{\Gamma} = (v_1)|_{\Gamma} - (v_2)|_{\Gamma}$, where $v_m = v|_{\Omega_m}$ is the restriction of $v$ to $\Omega_m$.
The first condition \eqref{eq2} is derived from the conservation of mass principle. The condition in \eqref{eq3} is called \emph{Henry condition}, c.f. \cite{Slattery,Sadhal, Bothe03VOF,Bothe04} and accounts for difference in the solubilities of the species in the different fluids. The coefficient $\beta = \beta(\bx,t)$ is piecewise constant:
\[ 
 \beta= \beta_m(t) >0 \quad \text{in}~~\Omega_m(t).
\]
As the solubilities do not coincide, in general we have $\beta_1 \neq \beta_2$ which leads to a jump discontinuity of $u$ across the interface. We assume that the ratio between $\beta_1$ and $\beta_2$ is bounded by a constant $c>0$:
$$
  \frac{1}{c} \beta_1 \leq \beta_2 \leq c \beta_1.
$$
The fact that the solution has a \emph{discontinuity across the (moving) interface} is the demanding aspect of the problem. 


To close the model we add initial and boundary conditions:
  \begin{eqnarray}
    u(\cdot , 0) & = &u_0 \quad \!\!\! \text{in}~~\Omega_m(0), ~m=1,2, \label{eq4} \\
    u(\cdot,t)& =&0  \quad \text{on}~~\partial \Omega, ~~t \in [0,T]. \label{eq5}
  \end{eqnarray}
\end{subequations}

For the treatment of the problem the case of a \emph{non-stationary interface} is significantly different from the case of a \emph{stationary interface}. In this paper we discuss several aspects of the finite element discretization method for the \emph{non-stationary interface} proposed in \cite{ReuskenLehrenfeld13}. 

The method in \cite{ReuskenLehrenfeld13} is based on a Nitsche-XFEM method. The Nitsche-XFEM method has been introduced and studied for different applications in \cite{Hansbo02,Hansbo03,Hansbo04,Hansbo05,Hansbo09,Hansbo11}. For such transport problems as in \eqref{eq1}-\eqref{eq5} the Nitsche-XFEM method has been investigated in \cite{Hansbo02,ReuskenHieu,ReuskenLehrenfeld12}. In \cite{Hansbo02} the Nitsche-XFEM method is analyzed for the simpler situation of a stationary problem without a discontinuity, i.e. $\beta_1=\beta_2$ and without convection, i.e. $\bw=0$. This method is extended to the instationary case with a stationary interface, a jump discontinuity across the interface and a convection velocity $\bw \neq 0$ in \cite{ReuskenHieu}. The problem of stabilizing the method in \cite{ReuskenHieu} w.r.t. a dominating convection term is discussed in \cite{ReuskenLehrenfeld12}.
In all these papers, the Nitsche-XFEM method is analyzed only for the case of a \emph{stationary interface}. In \cite{zunino13} a parabolic problem with a moving interface is discretized with a backward Euler/Nitsche-XFEM method and analyzed. This paper however only covers the case without discontinuities, i.e. $\beta_1=\beta_2$ and does not consider convection terms. 

In \cite[Chapter 10]{GReusken2011} it is shown that a well-posed weak formulation of the parabolic problem  \eqref{eq1}-\eqref{eq5} is most naturally formulated in terms of a space-time variational formulation. 
A combination of the Nitsche-XFEM method with a suitable \emph{space-time discontinuous Galerkin method}, 
the space-time XFEM-DG method is presented and analyzed in \cite{ReuskenLehrenfeld13}. 
The discretization has a (proven) second order accuracy in space and time for the transport problem in \eqref{eq1}-\eqref{eq5}.

An important characteristic of the Nitsche-XFEM method is the fact that the triangulation is not fitted to the interface. This also holds true for the space-time XFEM-DG method which renders it an Eulerian method.

The idea of combining space-time and extended finite elements has already been presented in \cite{Belytschko04} and applied to spatially one dimensional conservation equations without any analysis. 

In implementations of the spatially three dimensional case one has to deal with four dimensional geometries. In the literature space-time discretizations which explicitly construct four dimensional triangulations are rare. In \cite{Behr08} a decomposition of four dimensional prisms into simplices is applied to achieve real space-time adaptivity within one time slab. More global four dimensional simplex triangulations are considered in \cite{Neumueller11, NeumuellerPhD}.

The main new contributions of this paper are the following two:
\begin{enumerate}
\item A strategy to decompose the space-time domain which is intersected by an implicitly defined space-time interface into simplices and
\item numerical studies for the Nitsche XFEM-DG method in three space dimensions.
\end{enumerate}
As in  \cite{ReuskenLehrenfeld13}, we restrict to the diffusion dominated case, i.e., no stabilization w.r.t. convection is needed. 

\paragraph{structure and content of this paper}
In section~\ref{sectderivation} we present the discretization method that is considered in the analysis in \cite{ReuskenLehrenfeld13}.
This method is often not feasible in practice, due to the fact that it is assumed that volume integrals close to the space-time interface and surface integrals over the space-time interface are evaluated exactly. 
In practice a strategy for approximating these integrals will be necessary.
Implementation aspects are discussed in section~\ref{sectImpl}. Especially a strategy to construct a piecewise planar space-time interface approximation and corresponding polygonal approximations of the space-time subdomains is presented. A crucial point in this strategy is the decomposition of prisms (intersected by a piecewise planar interface) into (uncut) simplices. In the spatially three dimensional case this involves four dimensional geometries. A main contribution of this paper is a solution algorithm for this problem that is presented separately in section~\ref{sec:decomp4d}.
In section~\ref{sectExp} numerical experiments in three spatial dimensions are presented and discussed.


\section{The Nitsche XFEM-DG discretization} \label{sectderivation}
In this section we present the discretization method. We introduce some notation. The space-time cylinder $Q = \Omega \times (0,T] \subset \rr^{d+1}$ is partitioned into $N$ time slabs $Q^n = \Omega \times I_n$ with $I_n := (t_{n-1},t_n]$ and $0 = t_0 < t_1 < \ldots < t_N = T$. The length of each time interval $I_n$ is $\Delta t_n = t_n - t_{n-1}$. For simplicity we assume $\Delta t = \Delta t_n = \mathrm{const}$.  The discretization allows to solve for the discrete solution time slab by time slab which leads to a computational structure comparable to other (implicit) time-stepping schemes. In the following we will consider the problem on one time slab $Q^n$, we will therefore omit the subscript $n$. 

Within each time slab we assume that the triangulation of the spatial domain $\Omega$ is a shape regular decomposition into simplices $\mT = \{ T \}$. The corresponding characteristic mesh size is denoted by $h$. 
Corresponding to this triangulation of $\Omega$, we have a canonical triangulation of $Q$ into $d+1$-dimensional prisms. 
This triangulation is denoted by $\T^\ast = \{ \mQ^T \}$ where for each prism we have $\mQ^T = T \times I_n$ for a corresponding $d$-dimensional simplex $T$.
Note, that for different time slabs the triangulation is allowed to change. Further, the triangulation is \emph{not fitted} to the interface $\Gamma(t)$ (cf. figure \ref{fig:timeslabs}). 
We introduce the space-time finite element space used for both trial and test functions in the discretization. Based on $V_h$, the spatial finite element space of continuous, piecewise linear functions (with zero boundary values on $\partial \Omega$), the space-time finite element space of continuous piecewise \emph{bilinear} (linear in time, linear in space) functions is given by
\begin{equation} \label{defStandFE}
   W := \{ \, v: Q \to \Bbb{R} ~|~ v(\bx,t)= \phi_0(\bx) + t \phi_1(\bx),~~ \phi_0, \phi_1 \in V_h\,\}.
\end{equation}
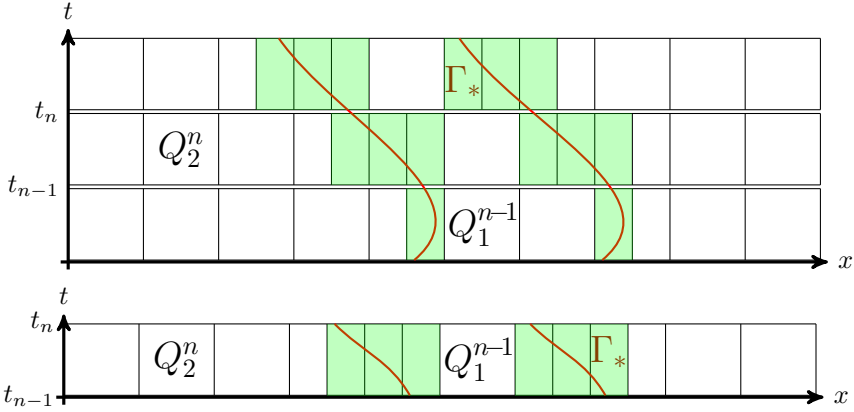
\begin{figure}
\centering
    \begin{tikzpicture}

      [

      rotate=0,

      ]

      \def\dw{0.025}


      \draw[->,very thick,>=stealth'] (-0.1,0)  -- (10.1,0) node(xline)[right]{$x$};

      \draw[->,very thick,>=stealth'] (0,-0.1)  -- (0,3.1) node(yline)[above]{$t$};


      \foreach \i in {0,...,10} {

        \foreach \j in {0,...,2} {

          \draw [] (\i,\j+\dw) -- (\i,1.0+\j-\dw); 

        }

      }

      \foreach \i in {0,...,2} {

        \draw [] (0,\i+\dw) -- (10,\i+\dw);

      }

      \foreach \i in {1,...,3} {

        \draw [] (0,\i-\dw) -- (10,\i-\dw);

      }
      \draw (0,2) node[left] {$t_n$};
      \draw (0,1) node[left] {$t_{n-1}$};

      \draw (5.625,2.4) node[left] {\Large \color{mildred} $\Gamma_\ast$};
      \draw (1.5,1.5) node[] {\Large $Q_2^n$};
      \draw (5.525,0.5) node[] {\Large $Q_1^{n\!-\!1}$};


      \draw [] (2.5,2+\dw) -- (2.5,3-\dw); 


      \draw [] (3.5,1+\dw) -- (3.5,2-\dw); 

      \draw [] (3.5,2+\dw) -- (3.5,3-\dw); 

      \draw [] (4.5,0+\dw) -- (4.5,1-\dw); 

      \draw [] (4.5,1+\dw) -- (4.5,2-\dw); 

      \draw [] (5.5,2+\dw) -- (5.5,3-\dw); 


      \draw [] (6.5,1+\dw) -- (6.5,2-\dw); 

      \draw [] (6.5,2+\dw) -- (6.5,3-\dw); 

      \draw [] (7.5,0+\dw) -- (7.5,1-\dw); 

      \draw [] (7.5,1+\dw) -- (7.5,2-\dw); 

      \draw[red,thick] (4.6,\dw) to[out=40,in=-55] (2.8,3-\dw);

      \draw[red,thick] (7.1,\dw) to[out=40,in=-55] (5.2,3-\dw);


      \fill[green,opacity=0.25] (4.5,\dw) rectangle (5,1-\dw); 

      \fill[green,opacity=0.25] (7,\dw) rectangle (7.5,1-\dw); 

      \fill[green,opacity=0.25] (3.5,1+\dw) rectangle (5,2-\dw); 

      \fill[green,opacity=0.25] (6,1+\dw) rectangle (7.5,2-\dw); 

      \fill[green,opacity=0.25] (2.5,2+\dw) rectangle (4,3-\dw); 

      \fill[green,opacity=0.25] (5,2+\dw) rectangle (6.5,3-\dw); 



    \end{tikzpicture}

    \begin{tikzpicture}

      [

      rotate=0,

      ]

      \def\dw{0.025}


      \draw[->,very thick,>=stealth'] (-0.1,0)  -- (10.1,0) node(xline)[right]{$x$};

      \draw[->,very thick,>=stealth'] (0,-0.1)  -- (0,1.1) node(yline)[above]{$t$};


      \foreach \i in {0,...,10} {

        \foreach \j in {0,...,0} {

          \draw [] (\i,\j+\dw) -- (\i,1.0+\j-\dw); 

        }

      }

      \foreach \i in {0,...,0} {

        \draw [] (0,\i+\dw) -- (10,\i+\dw);

      }

      \foreach \i in {1,...,1} {

        \draw [] (0,\i-\dw) -- (10,\i-\dw);

      }
      \draw (0,1) node[left] {$t_n$};
      \draw (0,0) node[left] {$t_{n-1}$};

      \draw (7.625,0.6) node[left] {\Large \color{mildred} $\Gamma_*$};
      \draw (1.5,0.5) node[] {\Large $Q_2^n$};
      \draw (5.525,0.5) node[] {\Large $Q_1^{n\!-\!1}$};


      \draw [] (3.5,0+\dw) -- (3.5,1-\dw); 

      \draw [] (4.5,0+\dw) -- (4.5,1-\dw); 

      \draw [] (6.5,0+\dw) -- (6.5,1-\dw); 

      \draw [] (7.5,0+\dw) -- (7.5,1-\dw); 

      \draw[red,thick] (4.6,\dw) to[out=120,in=-45] (3.6,1-\dw);

      \draw[red,thick] (7.2,\dw) to[out=120,in=-45] (6.2,1-\dw);


      \fill[green,opacity=0.25] (3.5,0+\dw) rectangle (5,1-\dw); 

      \fill[green,opacity=0.25] (6,0+\dw) rectangle (7.5,1-\dw); 

    \end{tikzpicture}
    \vspace*{-0.25cm}
    \caption{Sketch of the space-time domains $Q^n_m$. Note that within each time slab the triangulation has a tensor product structure $\T \times [t_{n-1},t_n]$. The triangulation is allowed to change between the time slabs.}
  \label{fig:timeslabs}
\end{figure}

The space-time interface $\Gamma_{\ast} := \bigcup_{t\in I_n} \Gamma(t)$ separates the time slab into its subdomains $Q_m := \bigcup_{t \in I_n} \Omega_m(t),\, m=1,2$. With $R_m$ the restriction operator on $L^2(Q)$, such that $R_m v = v|_{Q_m}$, we introduce the space-time XFEM space as 
\begin{equation} \label{XFEMdef}
  W^{\Gamma_\ast} := R_1 W \oplus R_2 W.
\end{equation}

To enforce the Henry condition \eqref{eq3} we use a space-time variant of the Nitsche method. To do this we need suitable jumps and averages across $\Gamma_{\ast}$. Within a prism-element $\mQ^T$, the averaging operator is defined as
\[
  \{v\}_{\Gamma_\ast}(\bx,t) := \kappa_1 ( v_1 )_{|\Gamma_\ast} + \kappa_2 ( v_2 )_{|\Gamma_\ast}, \quad (\bx,t) \in \Gamma_\ast 
\]
with $\mQ^T_m := Q_m \cap \mQ^T$ and the averaging weights
\[
  {\kappa_m}_{|\mQ^T}:= \frac{|\mQ^T_m|}{|\mQ^T|}.
\]
This is the space-time version of the averaging choice proposed in \cite{Hansbo02}.

We use a similar notation for the jump operators:
\[
 [v]_\Gs (\bx,t) = ( v_2 )_{|\Gs} - ( v_1 )_{|\Gs},\quad (\bx,t) \in \Gs.
\]
\begin{remark}
\rm
The choice for ${\kappa_m}|_{\mQ^T}$ is slightly different from the one in \cite{ReuskenLehrenfeld13}. As ${\kappa_m}|_{\mQ^T}$ is time-independent within one element $\mQ^T$ it is better suited for implementation. 
Under mild assumptions on the (approximated) space-time interface the theoretical results in \cite{ReuskenLehrenfeld13} carry over with only minor adjustments. 
\end{remark}

Further, for the weak enforcement of continuity in time (using the discontinuous Galerkin (DG) method) we use the definition $v_+^{n-1}(\cdot):= \lim_{\epsilon \downarrow 0} v(\cdot,t_{n-1}+\epsilon)$.

To transform the iterated integrals which appear in the Nitsche formulation to an integral on the space-time interface we use the following transformation formula:

\[
 \int_{t_{n-1}}^{t_n} \int_{\Gamma(t)} f(\bs,t) \, d\bs \, dt = \int_{\Gs} f(\bs) \big( 1+ (\bw \cdot \bn_\Gamma)^2 \big)^{-\frac12} \, d\bs
=: \int_\Gs f(\bs) \nu(\bs) \, d\bs,  
\]
 with $\nu(\bs)= \big( 1+ (\bw \cdot  \bn_\Gamma)^2 \big)^{-\frac12}$ and $\bw$ the (interface) velocity. Note that both the surface measure on $\Gamma(t)$ as well as on $\Gs$ are denoted with $d\bs$. 
Under the assumption that the space-time interface is sufficiently smooth, there holds for a constant $c_0>0$
\[
  c_0 \leq \nu(\bs) \leq 1 \quad \text{for all} ~~ \bs \in \Gs. 
\]
Below in the Nitsche bilinear form we use a weighting with $\nu(\bs)$.\\[1ex]
With the introduced notation we formulate the discrete variational formulation on the space-time slab $Q$. We define bilinear forms  $a(\cdot,\cdot)$, $b(\cdot,\cdot)$ and $N_{\Gs}(\cdot,\cdot)$ corresponding to the partial differential equation, the weak enforcement of continuity in time (DG) and the Nitsche terms to weakly impose the interface conditions, respectively. Further we define the linear forms $f(\cdot)$ and $c(\cdot;\cdot)$ corresponding to the right hand side source term $f$ and the initial conditions $u^{n-1}$:

%
\vspace*{-0.5cm}
\begin{eqnarray*}
  a(u,v) & = & \sum_{m=1}^2 \int_{Q_m} \big(\frac{\partial u_m}{\partial t} + \bw \cdot \nabla u_m\big) \beta_m v_m + \alpha_m \beta_m \nabla u_m \cdot \nabla v_m  \, d\bx, \\
  N_\Gs(u,v) & = &- \int_{\Gs} \nu(\bs) \{ \alpha \nabla u \cdot \bn\}_\Gs[\beta v]_\Gs \, d\bs - \int_{\Gs} \nu(\bs) \{ \alpha \nabla v \cdot \bn\}_\Gs[\beta u]_\Gs\, d\bs \\
  & &+ \lambda h_T^{-1} \int_{\Gs} \nu(\bs) [\beta u]_\Gs [\beta v]_\Gs \, d\bs,  \\
 b(u,v) &=& \sum_{m=1}^2 \int_{\Omega_{m}(t_{n-1})} \hspace*{-0.9cm} \beta_m u_+^{n-1} \,  v_+^{n-1} \, d\bx \\
  f(v) &=&  \sum_{m=1}^2 \int_{Q_m} \beta_m f \, v \, d\bx, \\
  c(\bar{u};v) &=& \sum_{m=1}^2 \int_{\Omega_{m}(t_{n-1})} \hspace*{-0.9cm} \beta_m u^{n-1} \, v_+^{n-1} \, d\bx.
\end{eqnarray*}
Note that both the volume measure on $\Omega_m(t)$ as well as on $Q_m$ are denoted with $d\bx$. 
Further, the (space-time) volume integrals are weighted with the factor $\beta$. 
In the definition of $c(\cdot,\cdot)$ the notation $\bar{u}$ should emphasize that $u^{n-1}$ is given (initial) data for the computations on the current time slab $Q$, e.g. for $n=1$, $u^{n-1}$ is the initial data $u^0$ of \eqref{eq4}. 
In $N_\Gs(\cdot,\cdot)$ the parameter $\lambda \geq 0$ has to be chosen sufficiently large. 

These bi- and linear forms are well-defined on the space-time XFEM space $W^{\Gamma_\ast}$ and the problem within one time slab reads as: Find $u \in W^{\Gamma_\ast}$, s.t.
\begin{equation} \label{oneslabprob}
 B(u,v) := a(u,v) + b(u,v) + N^\Gs(u,v) = f(v) + c(\bar{u},v) \qquad \forall \ v \in W^{\Gamma_\ast}
\end{equation}
Note that the solution of \eqref{oneslabprob} is the result of one time step within a \emph{time stepping} scheme.
 In \cite{ReuskenLehrenfeld13} it has been shown, that the presented method is a stable discretization scheme and a second order (in time and space) error bound for the $L^2(\Omega(T))$-norm has been proven. 
We recall the main result from \cite{ReuskenLehrenfeld13}.
\begin{theorem} \label{thm:mainresult}
Let $Q_m^{[0,T]} = \bigcup_{t\in[0,T]} \Omega_m(t)$, $m=1,2$ and $u \in H^2(Q_1^{[0,T]}\cup Q_2^{[0,T]})$ denote the solution of \eqref{eq1}-\eqref{eq5}. Further, let $U$ be the discrete solution of \eqref{oneslabprob} on every timeslab $Q^n$, $n=1,..,N$. Under certain assumptions on the regularity of the solution of the homogeneous backward problem to \eqref{eq1}-\eqref{eq5} (cf. \cite{ReuskenLehrenfeld13} for precise statement) the following holds:
$$
 \Vert (u - U)(\cdot,T) \Vert_{L^2(\Omega)} \leq c (h^2 + \Delta t^2) \Vert u \Vert_{H^2(Q_1^{[0,T]}\cup Q_2^{[0,T]})}.
$$
\end{theorem}
\begin{remark}[Linear systems] \label{rem:linsys}
\rm
Linear systems arising from a Nitsche-XFEM discretization may become ill-conditioned depending on the position of the interface. For the stationary case this has been discussed in \cite{burmanhansbo12,burmanzunino12}. This is also true for the discretization considered in this paper. However in our experience a remedy to this problem is a simple diagonal scaling. In the numerical examples in section \ref{sectExp} we solved the arising linear systems with a Jacobi-preconditioned GMRES method. In all these examples the iteration counts stayed within reasonable bounds. We do not know of any literature that rigorously discusses this issue without introducing additional stabilizations even for the stationary case. This is topic of ongoing research. 
\end{remark}

\begin{remark}[Computational efficiency] \label{rem:compeff}
\rm
A concern with the discretization method might be that the number of degrees of freedom is approximately twice as much as for a comparable first-order (in time) method. 
However the fact that we have a (proven) second order convergence in time justifies, in our opinion, the additional computational overhead. 
A crucial point for the discussion of the computational efficiency is the question of how to solve the arising linear systems efficiently. 
At this point we do not have a satisfying answer. 
In the implementation used for the numerical examples in section \ref{sectExp} we use standard iterative solvers and preconditioners (see also remark \ref{rem:linsys}) which are not adapted to the space-time context of the method.
The question of how to design more sophisticated solvers which exploit the space-time structure of the problem is the topic of ongoing research. A promising approach has been presented in \cite{bastingweller13} for the case without an interface. It is however not clear if these ideas can be applied in the context of XFEM methods. 
\end{remark}


\section{Implementation aspects}
\label{sectImpl}
In order to implement the method presented in section \ref{sectderivation} one has to compute the matrices and vectors representing the bi- and linear forms. In section \ref{sec:inttypes} we explain which types of integrals need to be calculated for that.

In practice the space-time interface $\Gs$ is typically not given explicitly, but implicitly, e.g. as the zero-level of a level-set function. In our application the level-set function is a piecewise quadratic function in space, and typically only given at discrete time levels $t_n$. 
In order to apply quadrature on the space-time objects, we want to have an explicit representation of the space-time interface. As this is practically not feasible, we construct an appropriate approximation $\Gsh$ which has an explicit representation. Such an approximation is discussed in section \ref{sec:stifapprox}. 
A different approach is presented in the recent paper \cite{MueKumObe13} where quadrature rules for \emph{implicitly} given domains by means of moment-fitting are derived. 

The special tensor product structure of each time slab is reflected in the construction and representation of finite element shape functions. A short discussion on the shape functions, especially w.r.t. XFEM can be found in section \ref{sec:trafofel}.

Once, the (space-time) geometries and finite element shape functions are defined one needs suitable quadrature rules for (regular) prisms and simplices. As these are discussed in many standard references (e.g. \cite{Stroud72}) for the case $d+1 \leq 3$, but only rarely in the $3+1$-dimensional case, section \ref{sec:onephaseint} and section \ref{sec:calcints} address this issue. In section \ref{sec:onephaseint} tensor product quadrature rules for prisms are applied for elements that are not intersected. The more involved situation when elements are cut by $\Gamma_{\ast,h}$ make use of the decomposition of an intersected $d+1$-prism into $d+1$-simplices. This is discussed in \ref{sec:calcints}.

\subsection{Integral types} 
\label{sec:inttypes}
Every bi- and linear form in \eqref{oneslabprob} has a natural decomposition into its element contributions, e.g. $a(u,v) = \sum_{\mQ^T \in \mathcal{T^\ast}} a^T(u,v) $. We consider the task of computing the element contributions of the (bi-)linear forms $a(\cdot,\cdot)$, $b(\cdot,\cdot)$, $N_\Gs(\cdot,\cdot)$, $f(\cdot)$ and $c(\bar{u};\cdot)$. 

As we need to calculate (approximations of) integrals of different kinds, we categorize these integrals before we discuss their numerical treatment. We distinguish those integrals in terms of the sets $S$ we are integrating on. The cases are denoted as \texttt{case (m,n,o)} where \texttt{m} is the dimension of $S$, \texttt{n} is the co-dimension of $S$ and $ \texttt{o} \in \{ \texttt{c}, \texttt{n} \}$ describes if the set $S$ is intersected by the space-time interface $\Gamma_\ast$ (\texttt{o=c}) or not (\texttt{o=n}).
We recall the notation for a prism $\mQ^T = T\times I_n$ . Accordingly we define $\mQ_m^T = \mQ^T \cap Q_m$. 
\subsubsection{$a(\cdot,\cdot), f(\cdot)$: $d\!+\!1$-dimensional measure, co-dimension 0}
Integrals appearing on each element for $a(\cdot, \cdot)$ are integrals on $d+1$-dimensional objects like
 $$ \int_{\mQ^T_m} f \, d\bx = \int_{t_{n-1}}^{t_n} \int_{T_m(t)} f \, d\bx \, , \, \, \, \, \text{with } T_m(t) = T \cap \Omega_m(t) $$ 
We distinguish two different situations: The prism $\mQ^T_m$ is not intersected by the (approximated) interface, i.e. the prism is completely in one phase and thus the volume to integrate on is the prism itself. We consider this as \texttt{case (d+1,0,n)} where numerical integration can exploit the tensor product structure. If on the other hand the prism $\mQ^T_m$ is intersected by the (approximated) interface, the geometry $\mQ^m_T$ is much more difficult to handle. In that case $d\!+\!1$-dimensional quadrature on subsimplices has to be applied. This is denoted by \texttt{case (d+1,0,c)}.
\subsubsection{$b(\cdot,\cdot), c(\hat{u},\cdot)$: $d$-dimensional measure, co-dimension 0}
The integrals in the element contributions of $b(\cdot,\cdot)$ and $c(\bar{u};\cdot)$ have the form 
$$
\int_{T_m(t_{n-1})} \!\!\!\! \!\!\!\! f \, d\bx 
$$
and thus are $d$-dimensional measures. Also here, we distinguish the case of a one phase element (i.e. an element which is not intersected), denoted by \texttt{case (d,0,n)} and the case of an intersected element, \texttt{case (d,0,c)}.
\subsubsection{$N^\Gs(\cdot, \cdot)$: $d$-dimensional measure, co-dimension 1}
For the space-time integrals stemming from the Nitsche stabilization bilinear form $N^\Gs(\cdot, \cdot)$ on each element we get terms like 
$$
 \int_{\Gamma_\ast^T} \nu(\bs) \cdot \, d\bs.
$$
where $\Gamma_\ast^T = \Gs \cap \mQ^T$.
Some terms also depend on the normal direction $\bn_\Gamma$. These integrals only appear on elements that are intersected. The measure is d-dimensional on the manifold $\Gs$ with co-dimension 1. This case is denoted as \texttt{case (d,1,c)}.
\subsection{Approximation of the space-time interface $\Gs$ and the space-time volumes $Q_m$}
\label{sec:stifapprox}
For the integration of space-time volumes and the space-time interface a discrete approximation of $\Gs$ and $Q_m$ which is feasible for numerical integration has to be found. This is relevant for the integration cases \texttt{(d+1,0,c)} and \texttt{(d,1,c)}. 

Note that case \texttt{(d,0,c)} is an integral on a d-dimensional simplex $T$. To deal with those we follow the strategies (for $d \leq 3$) discussed in \cite[Section 7.3]{GReusken2011}.
In order to get approximations $\Gamma_{\ast,h}$ and $Q_{m,h}$ to the space-time interface and volumes (\texttt{case (d+1,0,c)})  we proceed similarly as in lower ($d$) dimensions. 

\begin{figure}[ht]
\newcommand\prism[4]{
        \def\sca{1.0}
        \def\ax{0.7*\sca}
        \def\ay{-0.25*\sca}
        \def\bx{2*\sca}
        \def\by{1*\sca}
        \def\cx{1*\sca}
        \def\cy{2*\sca}

        \def\dt{1.7}

        \def\in{0}
        \def\out{250}

        \def\m{#1}
        \def\n{#2}
        \pgfmathsetmacro\nn{\n-1}
        \def\ldx{1.0/\n}

        \foreach \x in {0,...,\m} {
          \def\y{\x*\dt/\m}
          \foreach \r in {0,...,\nn} {
            \foreach \s in {0,...,\r} {
              \pgfmathsetmacro\aax{(\nn-\r)*\ldx}
              \pgfmathsetmacro\bbx{(\nn-\r)*\ldx+\ldx}
              \pgfmathsetmacro\ccx{(\nn-\r)*\ldx}
              \pgfmathsetmacro\aaz{\s*\ldx}
              \pgfmathsetmacro\bbz{\s*\ldx}
              \pgfmathsetmacro\ccz{\s*\ldx+\ldx}
              \pgfmathsetmacro\acoordx{\aax * \bx + \aaz  * \cx + (1.0-\aax-\aaz) * \ax}
              \pgfmathsetmacro\acoordz{\aax * \by + \aaz  * \cy + (1.0-\aax-\aaz) * \ay}
              \pgfmathsetmacro\bcoordx{\bbx * \bx + \bbz  * \cx + (1.0-\bbx-\bbz) * \ax}
              \pgfmathsetmacro\bcoordz{\bbx * \by + \bbz  * \cy + (1.0-\bbx-\bbz) * \ay}
              \pgfmathsetmacro\ccoordx{\ccx * \bx + \ccz  * \cx + (1.0-\ccx-\ccz) * \ax}
              \pgfmathsetmacro\ccoordz{\ccx * \by + \ccz  * \cy + (1.0-\ccx-\ccz) * \ay}
              \draw[#3] (\acoordx,\y,\acoordz) -- (\bcoordx,\y,\bcoordz)  -- (\ccoordx,\y,\ccoordz) -- cycle;
            }
          }
        }
        \pgfmathsetmacro\mmo{\m-1}
        \ifnum #4=1
        {
        \foreach \x in {0,...,\mmo} {
          \def\y{\x*\dt/\m}
          \def\yy{\x*\dt/\m+\dt/\m}
          \pgfmathsetmacro\aalpha{0.7-\x/\m*0.1}
          \pgfmathsetmacro\abeta{0.3+\x/\m*0.1}
          \pgfmathsetmacro\agamma{0.7-\x/\m*0.6}
          \pgfmathsetmacro\adelta{0.3+\x/\m*0.6}
          \pgfmathsetmacro\balpha{0.7-(\x+1)/\m*0.1}
          \pgfmathsetmacro\bbeta{0.3+(\x+1)/\m*0.1}
          \pgfmathsetmacro\bgamma{0.7-(\x+1)/\m*0.6}
          \pgfmathsetmacro\bdelta{0.3+(\x+1)/\m*0.6}

          \draw[red, thick, fill=red, opacity=0.15] 
             (\aalpha*\ax+\abeta*\bx,\y,\aalpha*\ay+\abeta*\by)
          to[out=\out,in=\in] (\agamma*\ax+\adelta*\cx,\y,\agamma*\ay+\adelta*\cy)
          -- (\bgamma*\ax+\bdelta*\cx,\yy,\bgamma*\ay+\bdelta*\cy)
          to[out=\in,in=\out] (\balpha*\ax+\bbeta*\bx,\yy,\balpha*\ay+\bbeta*\by)
          --cycle;

          \draw[red, thick, opacity=0.5] 
             (\aalpha*\ax+\abeta*\bx,\y,\aalpha*\ay+\abeta*\by)
          to[out=\out,in=\in] (\agamma*\ax+\adelta*\cx,\y,\agamma*\ay+\adelta*\cy);

          \draw[red, thick, opacity=0.5] 
            (\bgamma*\ax+\bdelta*\cx,\yy,\bgamma*\ay+\bdelta*\cy)
          to[out=\in,in=\out] (\balpha*\ax+\bbeta*\bx,\yy,\balpha*\ay+\bbeta*\by);

        }
        }
        \else
        {
        }
        \fi
        
        \foreach \x in {0,...,\mmo} {
          \def\yy{\x*\dt/\m}
          \pgfmathsetmacro\yyy{\yy+\dt/\m}

          \foreach \r in {0,...,\n} {
            \foreach \s in {0,...,\r} {
              \pgfmathsetmacro\xx{(\n-\r)*\ldx}
              \pgfmathsetmacro\zz{\s*\ldx}
              \pgfmathsetmacro\coordx{\xx * \bx + \zz  * \cx + (1.0-\xx-\zz) * \ax}
              \pgfmathsetmacro\coordz{\xx * \by + \zz  * \cy + (1.0-\xx-\zz) * \ay}
              \draw[#3] (\coordx,\yy,\coordz) -- (\coordx,\yyy,\coordz);
            }
          }
        }
        \foreach \x in {0,...,\m} {
          \def\yy{\x*\dt/\m}
          \pgfmathsetmacro\yyy{\yy+\dt/\m}
          \foreach \r in {0,...,\n} {
            \foreach \s in {0,...,\r} {
              \pgfmathsetmacro\xx{(\n-\r)*\ldx}
              \pgfmathsetmacro\zz{\s*\ldx}
              \pgfmathsetmacro\coordx{\xx * \bx + \zz  * \cx + (1.0-\xx-\zz) * \ax}
              \pgfmathsetmacro\coordz{\xx * \by + \zz  * \cy + (1.0-\xx-\zz) * \ay}
            }
          }
        }
}

\centering
\begin{tabular}{ccccc}
\hspace*{-0.5cm}
\begin{minipage}{0.3\textwidth}
\centering
\begin{tikzpicture}
		[scale=1.75,
                cube/.style={very thick,black,fill=green,opacity=0.25},
                incuboid/.style={very thick,black,fill=blue,opacity=0.25},
                grid/.style={very thin,gray},
                axis/.style={->,blue,thick}]
                \prism{1}{1}{thick}{1}
                \prism{1}{1}{thick}{0}

\end{tikzpicture}
\end{minipage}
\hspace*{-0.5cm}
&
$\rightarrow$
&
\hspace*{-0.5cm}
\begin{minipage}{0.3\textwidth}
\centering
\begin{tikzpicture}
		[scale=1.75,
                cube/.style={very thick,black,fill=green,opacity=0.25},
                incuboid/.style={very thick,black,fill=blue,opacity=0.25},
                grid/.style={very thin,gray},
                axis/.style={->,blue,thick}]
                \prism{2}{1}{densely dotted, thick}{1}
                \prism{1}{1}{thick}{0}
\end{tikzpicture}
\end{minipage}
\hspace*{-0.5cm}
&
$\rightarrow$
&
\hspace*{-0.5cm}
\begin{minipage}{0.3\textwidth}
\centering
\begin{tikzpicture}
		[scale=1.75,
                cube/.style={very thick,black,fill=green,opacity=0.25},
                incuboid/.style={very thick,black,fill=blue,opacity=0.25},
                grid/.style={very thin,gray},
                axis/.style={->,blue,thick}]
                \prism{2}{2}{densely dotted, thick}{1}
                \prism{1}{1}{thick}{0}
\end{tikzpicture}
\end{minipage}
\hspace*{-0.5cm}
\\
$m_t=1$, $m_s=1$ & & $m_t=2$, $m_s=1$ & & $m_t=2$, $m_s=2$ 
\end{tabular}
\vspace*{-0.2cm}
\caption{An intersected prism in $d+1$ dimensions, with $d=2$. The original prism (left), the prisms after uniform subdivision in time (middle) and after subdivision in space and time (right). }
\label{fig:subdiv}
\end{figure}
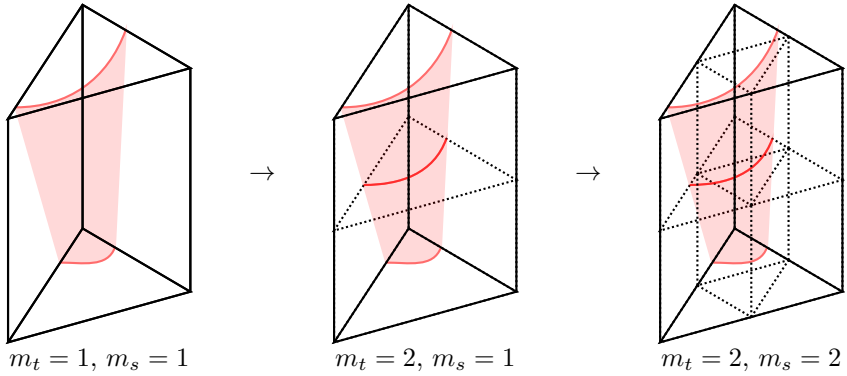

We consider the prism $ \mQ^T$ with a characteristic spatial length $h$ of $T$ (e.g. diameter) and the time step size $\Delta t = t_n - t_{n-1}$. We apply regular subdivisions in time and space. Each edge of $T$ is divided into $m_s$ parts of equal length and the time interval is divided into $m_t$ parts (see Figure \ref{fig:subdiv}). We get $m_t \cdot m_s^d$ smaller prisms $\{\mQ_i\}$ with spatial resolution $h/m_s$ and temporal resolution $\Delta t / m_t$.
%

%

Each (smaller) prism $\mQ_i$ is subdivided into $d+1$ $(d\!+\!1)$-simplices $\{ \mathcal{P}_j \}$ (cf. section \ref{sec:decomprefprism}). On $\mathcal{P}_j$ the level-set function is interpolated as a linear function in space-time (by simply evaluating the vertex values only).  As the level-set function is now represented as a linear function on each simplex, the according approximation of the zero-level of the level-set function is piecewise planar. 

\subsection{Finite elements}\label{sec:trafofel}
We briefly explain the construction of the basis functions for $W^{\Gamma_\ast}$. Let $\mathcal{J}^W = \{ (i,k), i\in\{ 1,..,N_V \}, k=1,2\}$ with $N_V$ the number of vertices of the spatial mesh denote the index set corresponding to basis functions $q_{i,k}$ in $W$. Based on $\mathcal{J}^W$ we can define the index subset of basis functions in $W$ which have to be ``enriched''. These are all basis functions which have a non-zero trace on the space-time interface:
\[
 \mathcal J_{\Gamma_\ast}:=\{\, (i,k) \in \mathcal J^W ~|~{\rm meas}_d\big(\Gamma_\ast\cap\supp (q_{i,k})\big)>0\}.
\]
For every index $(i,k) \in  \mathcal J_{\Gamma_\ast}$ we add the following basis function
\begin{equation} 
 q^{\Gamma_\ast}_{i,k} (\bx, t):= q_{i,k}(\bx, t)  \cdot (H_{\Gamma_\ast}(\bx,t)-H_{\Gamma_\ast}(\bx_i,t_{n-2+k})) , \quad (j, k) \in\mathcal J_{\Gamma_\ast}
\end{equation}
where $H_{\Gamma_\ast}$ is the characteristic function of $Q_2$ with $ H_{\Gamma_\ast}(\bx,t)=1 $ if $(\bx,t) \in Q_2$ and zero otherwise.
It is ensured that $q^{\Gamma_\ast}_{i,1}(\bx_i, t_{n-1})=0$, $q^{\Gamma_\ast}_{i,2}(\bx_i, t_{n})=0$  holds for all vertices $\bx_i$.
Using this construction it is obvious that the XFEM space $W^{\Gamma_\ast}$ can be decomposed into the standard space $W$ and additional basis functions which are discontinuous across the space-time interface:
\[
  W^{\Gamma_\ast} = W \oplus {\rm span} \big\{\, q_{i,k}^{\Gamma_\ast} ~|~ (i,k) \in \mathcal{J}_{\Gamma_\ast}~\big\}.
\]

Let the local index set of the local standard basis functions on one prism $\mQ^T$ be $\mathcal{J}^T=\{(i,k), i = 1,..,4, k = 1,2 \}$. The according shape functions are 
$$
q_{i,k}(\bx,t) := \phi_i(\bx) \psi_k(t), \quad (i,k) \in \mathcal{J}^T
$$
with $\phi_i(\bx) = \lambda_i$ where $\lambda_j$ denotes the barycentric coordinate corresponding to vertex $j$ inside of tetrahedron $T$ and $\psi_1(t) = \frac{1}{\Delta t} (t_n - t)$, $\psi_2(t) = \frac{1}{\Delta t} (t - t_{n-1})$.   

\subsection{ Numerical integration on non-intersected (space-time) volumes (case \texttt{(d+1,0,n)})}
\label{sec:onephaseint}
Whenever a (tetrahedral) element $T$ is not intersected by the interface $\Gamma(t)$ for the whole time interval $(t_{n-1},t_n]$ the volume integrals of $a(\cdot,\cdot)$ and $f(\cdot)$ act on (the complete) prismatic element $\mQ^T$. Consider for example the diffusion part on $\mQ^T$. One corresponding matrix entry for $I\!=\!(i,k),J\!=\!(j,l) \in \mathcal J^W$ is  $G^T_{I,J}:=g^T(q_{i,k},q_{j,l})$ with $g^T(u,v) := \int_{\mQ^T} \alpha_m \beta_m \nabla u \nabla v \, d\bx\, dt$ and can be computed using iterated integrals:
\begin{eqnarray*}
g^T(q_{i,k},q_{j,l}) 
& = & \int^{t_n}_{t_{n-1}} \!\! \int_T \alpha_m \beta_m (\nabla \phi_i) \psi_k (\nabla \phi_j) \psi_l\, d\bx \, dt \\
& = & \int^{t_n}_{t_{n-1}} \!\!\! \psi_k \psi_l dt \int_T \! \alpha_m \beta_m \nabla \phi_i \nabla \phi_j  d\bx 
=  \Delta t \ M_{k,l} \int_T \alpha_m \beta_m \nabla \phi_i \nabla \phi_j  d\bx
\end{eqnarray*}
where $\{ M_{k,l} \}_{k,l=1,2}$ is the one-dimensional mass matrix
$$
M = \left(\begin{array}{cc}
  1/3 & 1/6 \\ 1/6 & 1/3 
\end{array}\right).
$$
Thus, in this case the quadrature problem is reduced to a standard problem in $d$ dimensions. 
For the convection or r.h.s. term the velocity $\bw$ and the source term $f$ are not necessarily separable (in terms of $t$ and $\bx$). Nevertheless applying interpolation in time for $f$ and $\bw$, e.g.
$$f_{\Delta t} = f(t_{n-1},x) \cdot \psi_1(t) + f(t_{n},x) \cdot \psi_2(t) = \sum_{k=1}^2 f(t_k,x) \cdot \psi_k(t)  $$
or tensor product quadrature rules, e.g. 
\begin{eqnarray*}
h^T(q_{i,k},q_{j,l}) 
& := & \int^{t_n}_{t_{n-1}} \!\! \int_T \beta_m (\nabla \phi_i) \psi_k (\bw \cdot \nabla \phi_j) \psi_l\, d\bx \, dt \\
& \approx & \sum_{i=0}^{N_{ip}} \omega_i \psi_k(t_i) \psi_l(t_i) \int_T \beta_m (\nabla \phi_i) (\bw(t_i) \cdot \nabla \phi_j) \psi_l\, d\bx \, dt \\
\end{eqnarray*}
where $N_{ip}$, $t_i$ and $\omega_i$ are the information of the (1D) quadrature rule, reduces the complexity to a $d$-dimensional quadrature problem.

\subsection{Numerical integration on intersected space-time volumes and the space-time interface (case \texttt{(d+1,0,c)} and case \texttt{(d,1,c)})}
\label{sec:calcints}
If $\mQ^T$ is intersected some simplices $\mathcal{P}_j$ within the decomposition $\mQ^T = \{ \mP_j\} $ (see section \ref{sec:stifapprox}) are intersected by a planar approximation of the interface. Using the simplex and the (hyper-) plane one can find a decomposition of $\mathcal{P}_j$ into simplices $\{\mathcal{P}^{(k)}_j\}$ which are no longer intersected and form a decomposition of $\mathcal{P}_j$, $\mathcal{P}_j = \bigcup_k \mathcal{P}^{(k)}_j$. Furthermore the plane intersecting one simplex $\mathcal{P}_j$ can also be decomposed into uncut $d$-dimensional simplices. 
As this decomposition is neither obvious nor standard for the case $d=3$ a solution strategy is presented in detail in section \ref{sec:decomp4d}. 
Thus one can achieve an explicit decomposition of $\mQ^T_{m}$ into uncut $(d\!+\!1)$-dimensional simplices and of $\Gs_h$ into $d$-dimensional simplices. 
Once this decomposition is determined the integration can be applied simplex by simplex. 
Quadrature rules of high order for simplices can be found in standard references (see e.g. \cite{Stroud72}) if the dimension of the simplex is $d+1 \leq 3$.  For $d=3$, i.e. the simplex is four-dimensional (a pentatope) this is no longer standard. In the literature only a few integration rules can be found (see eg. \cite{Behr08} and \cite{Stroud72}). In appendix \ref{sec:quad4d} we quote lower order rules and strategies to generate higher order ones. Further in appendix \ref{sec:nu} we comment on the computation of the weighting factor $\nu(\bs), \, \bs \in \Gs$. 

\subsection{Numerical integration on intersected space volumes}
For the treatment of the cases \texttt{(d,0,n)} and \texttt{(d,0,c)} we refer to \cite[Section 7.3]{GReusken2011}.


\section{Numerical examples}
\label{sectExp}
In this section we present results of numerical experiments. Different from the experiments in \cite{ReuskenLehrenfeld13}, where only spatially one-dimensional situations have been considered, we consider spatially three-dimensional cases to assess the convergence behavior of the method. In this setting we restrict ourselves to piecewise linear (in time and space) finite element approximations.
In all examples we consider the $L^2(\Omega(T))$-error for different space and time resolutions. The time interval $[0,T]$ is divided into $n_t$ time slabs of equal size. Accordingly the time step size is $\Delta t= \frac{T}{n_t}$. The spatial domain is always a cuboid which is either a cube or divided into a small number of cubes. The cubes are divided into $(n_s)^d$ smaller cubes which are then divided into tetrahedra. The error behaviour is investigated w.r.t. refinements, i.e. series of $n_s$ and $n_t$.
For all computations we used a third order rule from \cite{Behr08} for the numerical integration on the (sub-)pentatopes. In appendix \ref{sec:quad4d} we review on integration rules for pentatopes. 
\subsection{Moving plane, quasi-1D}
\label{sec:tc5:onedim}
This case is the three-dimensional counterpart to the one-dimensional test case in \cite{ReuskenLehrenfeld13}.
The domain is the cube $\Omega=[0,2]^3$. The ``inner'' phase is contained in the domain $\Omega_1(t)=\{ \bx \in \Omega: |x_1-q(x_2,x_3)-r(t)| \leq D/2\}$, where $q: [0,2]^2 \rightarrow [0,2]$ is the graph describing the shape of the domain $\Omega_1$ and $r: [0,T] \rightarrow \rr$ the function describing the time-dependent shift of the interface in $x_1$-direction. $D = \frac23$ is the width of the domain in $x_1$-direction. The complementary domain is $\Omega_2(t)=\Omega \setminus \Omega_1(t)$. \\
The velocity field is given as $\bw = (\frac{\partial r}{\partial t}(t), 0, 0)$. As boundary conditions we apply periodicity in all directions, $u(x_i\!=\!0) = u(x_i\!=\!2), i=1,2,3$. This renders the problem essentially one-dimensional if $q(x_2,x_3)=const$. \\
We prescribe the r.h.s. source term $f$, such that the solution is given by
$$
u(\bx,t) = \sin (k \pi t) \cdot U^{m}(x_1-q(x_2,x_3)-r(t)), \qquad \bx \in \Omega^{m}(t), \quad m=1,2
$$
with
\begin{equation} \label{eq:tc5U12}
U^1(y) = ay+by^3 \quad \text{and} \quad U^2(y) = \sin(\pi y) 
\end{equation}
where $a$ and $b$ are chosen such that the interface conditions hold. 

The diffusivities are $(\alpha_1,\alpha_2)=(1,2)$ and the Henry weights $(\beta_1,\beta_2)=(1.5,1)$, resulting in $a \approx 1.02728$ and $b \approx 6.34294$. The problem is considered in the time interval $[0,T]$ with $T=1$.
\begin{figure}
  \centering
  \begin{minipage}{0.485\textwidth}
    \centering
    \begin{tikzpicture}
      [scale=1.25,
      cube/.style={very thick,black,fill=green,opacity=0.25},
      incuboid/.style={very thick,black,fill=blue,opacity=0.25},
      cen/.style={very thick,grey,fill=grey,opacity=.75},
      grid/.style={very thin,gray},
      axis/.style={->,blue,thick}]

      \coordinate (aaa) at (0,0,0);
      \coordinate (baa) at (2,0,0);
      \coordinate (aba) at (0,2,0);
      \coordinate (bba) at (2,2,0);
      \coordinate (aab) at (0,0,2);
      \coordinate (bab) at (2,0,2);
      \coordinate (abb) at (0,2,2);
      \coordinate (bbb) at (2,2,2);

      \coordinate (ccc) at (0.6666666,0,0);
      \coordinate (ccd) at (0.6666666,0,2);
      \coordinate (cdc) at (0.6666666,2,0);
      \coordinate (cdd) at (0.6666666,2,2);
      \coordinate (dcc) at (1.3333333,0,0);
      \coordinate (dcd) at (1.3333333,0,2);
      \coordinate (ddc) at (1.3333333,2,0);
      \coordinate (ddd) at (1.3333333,2,2);

      \coordinate (cen1) at (1,0,0);
      \coordinate (cen2) at (1,0,2);
      \coordinate (cen3) at (1,2,0);
      \coordinate (cen4) at (1,2,2);

      \draw[->,thick, dotted] (4.0/3.0,0.3,1.7) -- (2,0.3,1.7) node[above]{$\bw$};;
      \draw[->,color=red] (0,1,1) -- (2,1,1) node[above]{$x_1$};;
      \draw[->,thick,color=purple!30!blue] (1,2.25,0)node[right]{$q$} -- (1,2,0) ;;

      \draw[cube] (aaa) -- (baa) -- (bba) -- (aba) -- cycle; 
      \draw[cube] (aab) -- (bab) -- (bbb) -- (abb) -- cycle; 
      \draw[cube] (aaa) -- (aab) -- (abb) -- (aba) -- cycle; 
      \draw[cube] (baa) -- (bab) -- (bbb) -- (bba) -- cycle; 
      \draw[cube] (aaa) -- (aab) -- (bab) -- (baa) -- cycle; 
      \draw[cube] (aba) -- (abb) -- (bbb) -- (bba) -- cycle; 

      \draw[cen]  (cen1) -- (cen3) -- (cen4) -- (cen2) -- cycle; 
      \draw[incuboid]  (ccc) -- (dcc) -- (ddc) -- (cdc) -- cycle; 
      \draw[incuboid]  (ccd) -- (dcd) -- (ddd) -- (cdd) -- cycle; 
      \draw[incuboid]  (ccc) -- (ccd) -- (cdd) -- (cdc) -- cycle; 
      \draw[incuboid]  (dcc) -- (dcd) -- (ddd) -- (ddc) -- cycle; 
      \draw[incuboid]  (ccc) -- (ccd) -- (dcd) -- (dcc) -- cycle; 
      \draw[incuboid]  (cdc) -- (cdd) -- (ddd) -- (ddc) -- cycle; 

      \draw[cube] (0,0,0) -- (0,0,2);
      \draw[cube] (0,2,0) -- (0,2,2);
      \draw[cube] (2,0,0) -- (2,0,2);
      \draw[cube] (2,2,0) -- (2,2,2);

    \end{tikzpicture}
  \end{minipage}
  \begin{minipage}{0.45\textwidth}
    \centering
    \begin{tikzpicture}[scale=1.75]
      \draw[thick,->, color=red] (-1.2,0) -- (1.2,0) node[right] {$y$}; 
      \draw[->] (0,-0.8) -- (0,0.9) node[above] {$U^1 / U^2$};
      \draw[dashed] (-1.0/3.0,0.95) -- (-1.0/3.0,-0.95) node[right] {$\Gamma$};
      \draw[dashed] (1.0/3.0,0.95) -- (1.0/3.0,-0.95) node[right] {$\Gamma$};

      \draw[color=green!50!black, thick, samples=150,domain=-1:-1.0/3.0] plot (\x,{sin(pi * \x r)})   node[right] {};
      \draw[color=blue, thick, samples=150,domain=-1.0/3.0:1.0/3.0] plot (\x,{1.02728 * \x + 6.34294*( \x * \x * \x)})   node[right] {};
      \draw[color=green!50!black, thick, samples=150,domain=1.0/3.0:1.0] plot (\x,{sin(pi * \x r)})   node[right] {};
    \end{tikzpicture} 
  \end{minipage}
  \vspace*{-0.2cm}
  \caption{Sketch of the geometrical setup of the numerical example in section \ref{sec:tc5:onedim} (left) , $\Omega_1$ is blue and $\Omega_2$ is green. The time-independent part of the solution $U^1 / U^2$ (see \eqref{eq:tc5U12}) is also sketched (right).}
\label{fig:tc5sketch}
\vspace*{-0.2cm}
\end{figure}
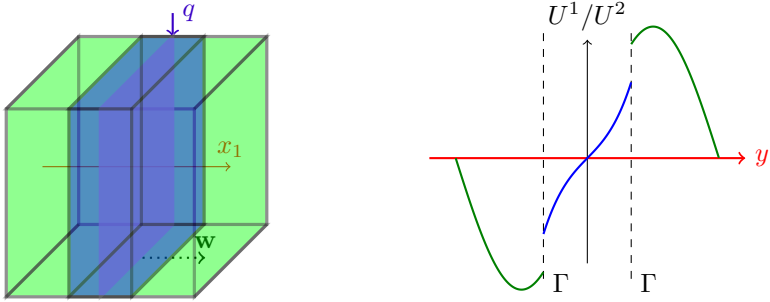
\subsubsection{Planar (in space and time) interface} \label{sec:tc5:onedim:linmove}
We choose $q(x_2,x_3)=1$ and $r(t)=0.25 t$, hence the space-time interface is planar. Thus the proposed method for the approximation of the space-time interface is exact for every $m_t, m_s \geq 1$, where $m_t$ and $m_s$ are the number of subdivisions in time and each space direction, respectively, cf. section \ref{sec:stifapprox}. We choose $m_s=m_t=1$.
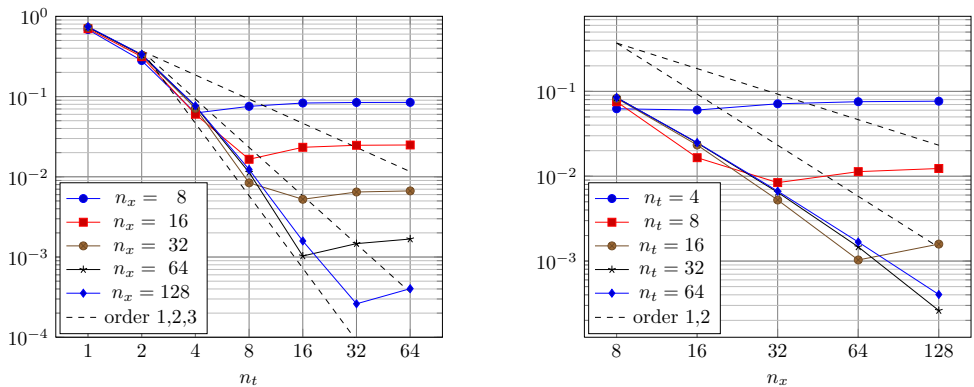
\begin{figure}[ht]
  \begin{tikzpicture}[scale=0.75]
    \begin{loglogaxis}[
      xlabel=$n_t$,
      ymin=1e-4,
      ymax=1,
      log basis x=2,
      grid=both,
      major grid style={black!50},
      legend style={
        at={(0.01,0.01)},
        anchor = south west
      },
      xticklabel=\pgfmathparse{2^\tick}\pgfmathprintnumber{\pgfmathresult}
      ]

      \addplot plot coordinates {
        ( 1 ,	0.682559 )
        ( 2 ,	0.28091 )
        ( 4 ,	0.0623801 )
        ( 8 ,	0.0754833 )
        ( 16 ,	0.0830545 )
        ( 32 ,	0.0843712 )
        ( 64 ,	0.0845707 )
      };
      \addplot plot coordinates {
        ( 1 ,	0.700457 )
        ( 2 ,	0.31375 )
        ( 4 ,	0.0602584 )
        ( 8 ,	0.0165317 )
        ( 16 ,	0.0233162 )
        ( 32 ,	0.0247002 )
        ( 64 ,	0.0249133 )
      };
      \addplot plot coordinates {
        ( 1 ,	0.727595 )
        ( 2 ,	0.327709 )
        ( 4 ,	0.0713872 )
        ( 8 ,	0.00840438 )
        ( 16 ,	0.00524292 )
        ( 32 ,	0.00646312 )
        ( 64 ,	0.0066706 )
      };
      \addplot plot coordinates {
        ( 1 ,	0.73301 )
        ( 2 ,	0.332598 )
        ( 4 ,	0.0755701 )
        ( 8 ,	0.0113221 )
        ( 16 ,	0.00103009 )
        ( 32 ,	0.00146948 )
        ( 64 ,	0.00167639 )
      };
      \addplot plot coordinates {
        ( 1 ,	0.742317 )
        ( 2 ,	0.334165 )
        ( 4 ,	0.0766494 )
        ( 8 ,	0.0123275 )
        ( 16 ,	0.00158459 )
        ( 32 ,	0.000261255 )
        ( 64 ,	0.000403262 )
      };
      \addplot plot[no marks,color=black,style=dashed] coordinates {
        ( 2 ,	0.3711585 )
        ( 4 ,	0.18557925 )
        ( 8 ,	0.092789625 )
        ( 16 ,	0.0463948125 )
        ( 32 ,	0.02319740625 )
        ( 64 ,	0.011598703125 )
      };
      \addplot plot[no marks,color=black,style=dashed] coordinates {
        ( 2 ,	0.3711585 )
        ( 4 ,	0.092789625 )
        ( 8 ,	0.02319740625 )
        ( 16 ,	0.0057993515625 )
        ( 32 ,	0.00144983789063 )
        ( 64 ,	0.000362459472656 )
      };
      \addplot plot[no marks,color=black,style=dashed] coordinates {
        ( 2 ,	0.3711585 )
        ( 4 ,	0.0463948125 )
        ( 8 ,	0.0057993515625 )
        ( 16 ,	0.000724918945313 )
        ( 32 ,	9.06148681641e-05 )
      };
      \legend{
        $n_x= \ \ 8 $\\
        $n_x= \ 16 $\\
        $n_x= \ 32 $\\
        $n_x= \ 64 $\\
        $n_x= 128 $\\
        order $1$,$2$,$3$\\}
    \end{loglogaxis}
  \end{tikzpicture}
  \hfill
  \begin{tikzpicture}[scale=0.75]
    \begin{loglogaxis}[
      xlabel=$n_x$,
      log basis x=2,
      grid=both,
      major grid style={black!50},
      legend style={
        at={(0.01,0.01)},
        anchor = south west
      },
      xticklabel=\pgfmathparse{2^\tick}\pgfmathprintnumber{\pgfmathresult}
      ]

      \addplot plot coordinates {
        ( 8 ,	0.0623801 )
        ( 16 ,	0.0602584 )
        ( 32 ,	0.0713872 )
        ( 64 ,	0.0755701 )
        ( 128 ,	0.0766494 )
      };
      \addplot plot coordinates {
        ( 8 ,	0.0754833 )
        ( 16 ,	0.0165317 )
        ( 32 ,	0.00840438 )
        ( 64 ,	0.0113221 )
        ( 128 ,	0.0123275 )
      };
      \addplot plot coordinates {
        ( 8 ,	0.0830545 )
        ( 16 ,	0.0233162 )
        ( 32 ,	0.00524292 )
        ( 64 ,	0.00103009 )
        ( 128 ,	0.00158459 )
      };
      \addplot plot coordinates {
        ( 8 ,	0.0843712 )
        ( 16 ,	0.0247002 )
        ( 32 ,	0.00646312 )
        ( 64 ,	0.00146948 )
        ( 128 ,	0.000261255 )
      };
      \addplot plot coordinates {
        ( 8 ,	0.0845707 )
        ( 16 ,	0.0249133 )
        ( 32 ,	0.0066706 )
        ( 64 ,	0.00167639 )
        ( 128 ,	0.000403262 )
      };
      \addplot plot[no marks,color=black,style=dashed] coordinates {
        ( 8 ,	0.3711585 )
        ( 16 ,	0.18557925 )
        ( 32 ,	0.092789625 )
        ( 64 ,	0.0463948125 )
        ( 128 ,	0.02319740625 )
      };
      \addplot plot[no marks,color=black,style=dashed] coordinates {
        ( 8 , 0.3711585 )
        ( 16 ,	0.092789625 )
        ( 32 ,	0.02319740625 )
        ( 64 ,	0.0057993515625 )
        ( 128 ,	0.00144983789063 )
      };
      \legend{
        $n_t= 4 $\\
        $n_t= 8 $\\
        $n_t= 16 $\\
        $n_t= 32 $\\
        $n_t= 64 $\\
        order $1$,$2$\\}
    \end{loglogaxis}
  \end{tikzpicture}
  \vspace*{-0.35cm}
  \caption{Convergence in $L^2(\Omega(T))$-norm w.r.t. refinements in time (left) and space (right) for test case in section \ref{sec:tc5:onedim:linmove}.}
  \label{fig:tc5:onedim:linmove}
\end{figure}

\begin{table}[ht]
  \small
  \begin{center}
    \begin{tabular}{|r|r|r|r|r||r|r|}
      \hline
      $\nt$ $\backslash$ $\ns$ & 8 & 16 & 32 & 64 & \ccg     128 & \ccg \eoct \\
      \hline
      1 & 0.6826 & 0.7005 & 0.72760 & 0.73301 & \ccg 0.742317 & \ccg \\
      2 & 0.2809 & 0.3138 & 0.32771 & 0.33260 & \ccg 0.334165 & \ccg 1.15 \\
      4 & 0.0624 & 0.0603 & 0.07139 & 0.07557 & \ccg 0.076649 & \ccg 2.12 \\
      8 & 0.0755 & 0.0165 & 0.00840 & 0.01132 & \ccg 0.012328 & \ccg 2.64 \\
      16 & 0.0831 & 0.0233 & 0.00524 & 0.00103 & \ccg 0.001585 & \ccg 2.96 \\
      32 & 0.0844 & 0.0247 & 0.00646 & 0.00147 & \ccg 0.000261 & \ccg 2.6 \\
      \hline
      \ccg 64 & \ccg 0.0846 & \ccg 0.0249 & \ccg 0.00667 & \ccg 0.00168 & \ccg 0.000403 & \ccg -0.63 \\
      \hline
      \ccg \eocs & \ccg & \ccg 1.76 & \ccg 1.90 & \ccg 1.99 & \ccg 2.06 & \ccg \\
      \hline
    \end{tabular}
 \end{center}
 \caption{$\Vert u_h-u \Vert_{L^2(\Omega(T))}$ for different refinements in time and space for the test case in section \ref{sec:tc5:onedim:linmove}. The last column shows the estimated order of convergence w.r.t. time (\eoct) using the finest spatial resolution ($\ns=128$), the last row shows the estimated order of convergence w.r.t. space (\eocs) using the finest temporal resolution ($\nt=64$).}
 \label{tbl:tc5:onedim:linmove}
\vspace*{-0.5cm}
\end{table}
%

%
In Table \ref{tbl:tc5:onedim:linmove} and Figure \ref{fig:tc5:onedim:linmove} we give the resulting error $\Vert u_h-u \Vert_{L^2(\Omega(T))}$ for different mesh and time step sizes. Similar to the results in \cite{ReuskenLehrenfeld13} we observe a third order convergence w.r.t. the time step size $\Delta t$ and a second order convergence w.r.t. the mesh size.
\begin{remark} \rm
  In \cite[Theorem 12.7]{Thomee97} for the corresponding DG-FEM method applied to the standard heat equation an error bound with third order convergence w.r.t. $\Delta t$ has been derived. Note however that the analysis does not carry over for the case of the Nitsche-XFEM discretization. 
\end{remark}
\subsubsection{Nonlinear moving interface} \label{sec:tc5:onedim:sinmovecurved}
We consider $q(x_2,x_3)=\frac78+\frac14 x_2^2 (2-x_2)^2$ and $r(t)=\frac{1}{4 \pi} \sin(2 \pi t)$ which leads to a space-time interface which is no longer planar. 
The geometrical approximation of the space-time interface in this paper is piecewise planar, i.e. the maximum distance between $\Gamma_\ast$ and its approximation $\Gamma_{\ast,h}$ converges with second order w.r.t. increasing $n_t \cdot m_t, n_s \cdot m_s$.
In Table \ref{tbl:tc5:onedim:sinmovecurved} the error $\Vert u_h-u \Vert_{L^2(\Omega(T))}$ on a fixed (fine) spatial grid with resolution $64\!\times\! 64\!\times\! 64$ for different numbers of time steps is listed. In order to investigate the impact of the approximation of $\Gamma_\ast$ we performed the computation with different numbers of subdivisions $m_s$, $m_t$. 
\begin{table}[ht]
  \small
  \begin{center}
    \begin{tabular}{|r|r|r|r|r|r|r|r|}
      \hline
      $\nt$ & 1 & 2 & 4 & 8 & \ccg 16 & 32 & 64 \\
      \hline
      $\ms=1,\mt=1$ & 2.50 & 2.89 & 0.547 & 0.137 & \ccg 0.0342 & 0.00879 & 0.00317 \\
      \hline
      \eoct &  & -0.21 & 2.40 & 2.00 & \ccg 2.00 & 1.96 & 1.49 \\
      \hline
      $\ms=1,\mt=2$ & 2.49 & 0.817 & 0.168 & 0.0374 & \ccg 0.00878 & 0.00301 & 0.00241 \\
      \hline
      \eoct &  & 1.61 & 2.28 & 2.17 & \ccg 2.09 & 1.54 & 0.32 \\
      \hline
      $\ms=1,\mt=4$ & 0.520 & 0.481 & 0.0985 & 0.0167 & \ccg 0.00284 & 0.00219 & 0.00236 \\
      \hline
      \eoct &  & 0.11 & 2.29 & 2.56 & \ccg 2.56 & 0.37 & -0.11 \\
      \hline
      $\ms=1,\mt=8$ & 0.491 & 0.412 & 0.0910 & 0.0143 & \ccg 0.00189 & 0.00212 & 0.00236 \\
      \hline
      \eoct &  & 0.25 & \ccg 2.18 & \ccg 2.67 & \ccg 2.92 & -0.17 & -0.15 \\
      \hline
      $\ms=4,\mt=8$ & 0.491 & 0.412 & 0.0909 & 0.0142 & \ccg 0.00179 & 0.00207 &  \\
      \hline
      \eoct &  & 0.25 & 2.18 & 2.68 & \ccg 2.99 & -0.20 &  \\
      \hline
    \end{tabular}
  \end{center}
  \caption{Error $\Vert u_h-u \Vert_{L^2(\Omega(T))}$  for different temporal refinements and quadrature subdivisions on a regular $64\!\times\! 64\!\times\! 64$ tetrahedral mesh for the test case in section \ref{sec:tc5:onedim:sinmovecurved}.}
  \label{tbl:tc5:onedim:sinmovecurved}
  \vspace*{-0.6cm}
\end{table}
%
%
%
The results, shown in Table \ref{tbl:tc5:onedim:sinmovecurved}, indicate an error bound behaviour of the form
$$
\Vert u_h - u \Vert_{L^2(\Omega(T))} \leq C_1 \Delta t^3 + C_2 \left(\Delta t / m_t \right)^2 + C_3(h) 
$$
where $C_1$ is independent of the approximation of $\Gamma_\ast$. $C_2$ is directly related to the interface approximation errors. If the interface approximation is exact (as in the last section) $C_2$ is zero. The function $C_3(h)$ describes the spatial error due to the method \emph{and} the piecewise linear interface approximation for the numerical integration. It is thereby the part of the error that can not be reduced by refinements in time. In this examples $C_3(h) \approx 0.002$. Furthermore, we observe that in this example $C_3$ is essentially independent of $m_s$.

For $m_t$ sufficiently large, i.e. $m_t  > \sqrt{C_2 / (C_1 \Delta t)}$ and $h$ sufficiently small, the first term dominates the error, such that one observes a third order in time convergence. This does not hold if $m_t$ is too small.
Especially for $m_s\!=\!m_t\!=\!1$, the error is converging with (only) second order, due to a dominating interface approximation error.
\iftoggle{extended}{%
  \begin{remark}\label{rem:xfemdofs}
    \rm
    In order to investigate the additional effort within one time step due to additional XFEM unknowns, we consider the ratio between the maximum number of extended (XFEM) unknowns and standard (space-time) finite element unknowns within one time slab. In Figure \ref{fig:tc5_xfem} a sketch of the corresponding situation is shown. 
    If the interface is well resolved, the number of unknowns close to the interface increases by a factor of $2^{d-1}$ for one uniform (spatial) refinement whereas the overall number of unknowns increases with $2^d$. Thus the ratio decreases linearly with the spatial resolution. In Table \ref{tbl:tc5:onedim:sinmovecurved:unkn} the corresponding numbers for this test case are given which are in agreement with the expected behaviour. 
  \end{remark}
  \input{figtabs/sec4_xfemref}
  \begin{table}[ht]
  \small
  \begin{center}
    \begin{tabular}{|r|r|r|r|r|r|r|}
      \hline
      Std. unkn. & 1024 & 8192 & 65536 & 524288 & 4194304 \\
      \hline
      $\nt$ $\backslash$ $\ns$ & 8 & 16 & 32 & 64 & 128 \\
      \hline
      1 & 736 (72\pr) & 3648  (45\pr) & 19328 (30\pr) & 119808 (23\pr) & 813056(19\pr) \\
      2 & 656 (64\pr) & 3008  (37\pr) & 14336 (22\pr) & 78848 (15\pr) & 479232(11\pr) \\
      4 & 656 (64\pr) & 3008  (37\pr) & 14336 (22\pr) & 78848 (15\pr) & 479232(11\pr) \\
      8 & 656 (64\pr) & 2880  (35\pr) & 13056 (20\pr) & 67072 (13\pr) & 383488 (9\pr) \\
      16 & 640 (63\pr) & 2624  (32\pr) & 11136 (17\pr) & 52992 (10\pr) & 275456 (7\pr) \\
      32 & 640 (63\pr) & 2496  (30\pr) & 10368 (16\pr) & 45824  (9\pr) & 198656 (5\pr) \\
      64 & 640 (60\pr) & 2496  (30\pr) & 9856 (15\pr) & 41472  (8\pr) & 182784 (4\pr) \\
      \hline
    \end{tabular}
  \end{center}
  \caption{Number of standard (space-time) unknowns (first row) and maximal number of additional XFEM unknowns for one time slab for different spatial and temporal resolutions for test case in section \ref{sec:tc5:onedim:sinmovecurved}. In brackets the ratio between XFEM and standard unknowns is added.}
  \label{tbl:tc5:onedim:sinmovecurved:unkn}
\end{table}


  %
}{%
}
\begin{remark}
  \rm
  To decrease the (space-time) interface approximation error one can either choose smaller time steps or a larger refinement factor $m_t$ for the construction of $\Gamma_{\ast,h}$. The computation with a fixed $\Delta t= \Delta \tilde{t}$ and $m_t=\tilde{m}_t>1$ is cheaper than a computation with $\Delta t = \Delta \tilde{t} / \tilde{m}_t$ and $m_t=1$. In Figure \ref{fig:temporalref} a sketch of both strategies is shown. For $m_t>1$  additional effort due to the decomposition strategy and quadrature within one time step is required. However if the interface is resolved, this is only required for a small number of elements\iftoggle{extended}{ (cf. Remark \ref{rem:xfemdofs})}{}. The number of time steps and thereby the number of linear systems that have to be solved however is reduced by a factor of $\tilde{m}_t$. Note that the solution of linear systems is typically the most time consuming part.
\end{remark}
\input{figtabs/sec4_subdivtime}
\subsection{Rotational symmetric solution on a moving sphere}
\label{sec:tc1:rot}
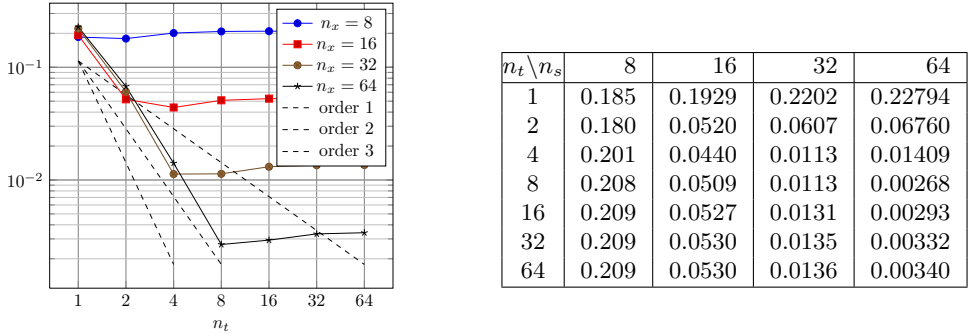
\begin{figure}[ht]
  \begin{minipage}{0.475\textwidth}
    \begin{tikzpicture}[scale=0.666]
      \begin{loglogaxis}[
        xlabel=$n_t$,
        log basis x=2,
        grid=both,
        major grid style={black!50},
        xticklabel=\pgfmathparse{2^\tick}\pgfmathprintnumber{\pgfmathresult}
        ]

        \addplot plot coordinates {
          ( 1 ,	0.185408 )
          ( 2 ,	0.179611 )
          ( 4 ,	0.201167 )
          ( 8 ,	0.207848 )
          ( 16 ,	0.209004 )
          ( 32 ,	0.209086 )
          ( 64 ,	0.209031 )
        };
        \addplot plot coordinates {
          ( 1 ,	0.192864 )
          ( 2 ,	0.0519731 )
          ( 4 ,	0.0439745 )
          ( 8 ,	0.0509014 )
          ( 16 ,	0.0526674 )
          ( 32 ,	0.0529913 )
          ( 64 ,	0.0530154 )
        };
        \addplot plot coordinates {
          ( 1 ,	0.220187 )
          ( 2 ,	0.0607429 )
          ( 4 ,	0.0112854 )
          ( 8 ,	0.0113351 )
          ( 16 ,	0.0131086 )
          ( 32 ,	0.0135095 )
          ( 64 ,	0.0135703 )
        };
        \addplot plot coordinates {
          ( 1 ,	0.227944 )
          ( 2 ,	0.0676023 )
          ( 4 ,	0.014087 )
          ( 8 ,	0.00267693 )
          ( 16 ,	0.00292528 )
          ( 32 ,	0.00332403 )
          ( 64 ,	0.00340398 )
        };
        \addplot plot[no marks,color=black,style=dashed] coordinates {
          ( 1 ,	0.113972 )
          ( 2 ,	0.056986 )
          ( 4 ,	0.028493 )
          ( 8 ,	0.0142465 )
          ( 16 ,	0.00712325 )
          ( 32 ,	0.003561625 )
          ( 64 ,	0.0017808125 )
        };
        \addplot plot[no marks,color=black,style=dashed] coordinates {
          ( 1 ,	0.113972 )
          ( 2 ,	0.028493 )
          ( 4 ,	0.00712325 )
          ( 8 ,	0.0017808125 )
        };
        \addplot plot[no marks,color=black,style=dashed] coordinates {
          ( 1 ,	0.113972 )
          ( 2 ,	0.0142465 )
          ( 4 ,	0.0017808125 )
        };
        \legend{
          $n_x= 8 $\\
          $n_x= 16 $\\
          $n_x= 32 $\\
          $n_x= 64 $\\
          order $1$\\
          order $2$\\
          order $3$\\}
      \end{loglogaxis}
    \end{tikzpicture}
  \end{minipage}
  $\ \ $
  \begin{minipage}{0.45\textwidth}
    \small
    \vfill
    \begin{tabular}{|@{}c@{}|r|r|r|r|r|r|}
      \hline
      $\nt$$\backslash$$\ns$ & 8 & 16 & 32 & 64 \\
      \hline
      1 & 0.185 & 0.1929 & 0.2202 & 0.22794 \\
      2 & 0.180 & 0.0520 & 0.0607 & 0.06760 \\
      4 & 0.201 & 0.0440 & 0.0113 & 0.01409 \\
      8 & 0.208 & 0.0509 & 0.0113 & 0.00268 \\
      16 & 0.209 & 0.0527 & 0.0131 & 0.00293 \\
      32 & 0.209 & 0.0530 & 0.0135 & 0.00332 \\
      64 & 0.209 & 0.0530 & 0.0136 & 0.00340 \\
      \hline
    \end{tabular}
    %
  \end{minipage}
  \vspace*{-0.2cm}
  \caption{Convergence in $L^2(\Omega(T))$-norm w.r.t. refinements in time and space as a plot and as tabulated values for test case in section \ref{sec:tc1:rot}.}\label{fig:tc1_conv}
\end{figure}

In this example we consider a rotational symmetric solution for a stationary sphere. The sphere is then translated with a time-dependent velocity. The time interval is $[0,T]$ with $T=0.5$ and the domain is the cube $\Omega=[0,2]^3$. 
\begin{figure}[ht]
\begin{minipage}{0.485\textwidth}
\centering
\begin{tikzpicture}
		[scale=1.25,
                cube/.style={very thick,black,fill=green,opacity=0.25},
                incuboid/.style={very thick,black,fill=blue,opacity=0.25},
                grid/.style={very thin,gray},
                axis/.style={->,blue,thick}]

        \coordinate (aaa) at (0,0,0);
        \coordinate (baa) at (2,0,0);
        \coordinate (aba) at (0,2,0);
        \coordinate (bba) at (2,2,0);
        \coordinate (aab) at (0,0,2);
        \coordinate (bab) at (2,0,2);
        \coordinate (abb) at (0,2,2);
        \coordinate (bbb) at (2,2,2);

        \coordinate (ccc) at (0.6666666,0,0);
        \coordinate (ccd) at (0.6666666,0,2);
        \coordinate (cdc) at (0.6666666,2,0);
        \coordinate (cdd) at (0.6666666,2,2);
        \coordinate (dcc) at (1.3333333,0,0);
        \coordinate (dcd) at (1.3333333,0,2);
        \coordinate (ddc) at (1.3333333,2,0);
        \coordinate (ddd) at (1.3333333,2,2);

	\draw[->,thick, dotted] (4.0/3.0,1,1) -- (2,1,1) node[above]{$\bw$};;
	\draw[->,color=red] (1,0,1) -- (1,2,1) node[above]{$y$};;

	\draw[cube] (aaa) -- (baa) -- (bba) -- (aba) -- cycle; 
	\shade[ball color=blue, opacity=0.5] (1,1,1) circle(0.45);
	\draw[cube] (aab) -- (bab) -- (bbb) -- (abb) -- cycle; 
	\draw[cube] (aaa) -- (aab) -- (abb) -- (aba) -- cycle; 
	\draw[cube] (baa) -- (bab) -- (bbb) -- (bba) -- cycle; 
	\draw[cube] (aaa) -- (aab) -- (bab) -- (baa) -- cycle; 
	\draw[cube] (aba) -- (abb) -- (bbb) -- (bba) -- cycle; 

	\draw[cube] (0,0,0) -- (0,0,2);
	\draw[cube] (0,2,0) -- (0,2,2);
	\draw[cube] (2,0,0) -- (2,0,2);
	\draw[cube] (2,2,0) -- (2,2,2);

\end{tikzpicture}
\end{minipage}
\hfill
\begin{minipage}{0.45\textwidth}
\centering
\begin{tikzpicture}[scale=1.75]
\draw[thick,->, color=red] (-1.2,0) -- (1.2,0) node[right] {$y$}; 
\draw[->] (0,-0.8) -- (0,1.1) node[above] {$U^1 / U^2$};
\draw[dashed] (-1.0/3.0,0.95) -- (-1.0/3.0,-0.75) node[above right] {$\Gamma$};
\draw[dashed] (1.0/3.0,0.95) -- (1.0/3.0,-0.75) node[above right] {$\Gamma$};

\draw[color=green!50!black, thick, samples=150,domain=-1:-1.0/3.0] plot (\x,{0.7*cos(pi * \x r)})   node[right] {};
\draw[color=blue, thick, samples=150,domain=-1.0/3.0:1.0/3.0] plot (\x,{0.7*1.1569 - 0.7*8.1621*( \x * \x)})   node[right] {};
\draw[color=green!50!black, thick, samples=150,domain=1.0/3.0:1.0] plot (\x,{0.7*cos(pi * \x r)})   node[right] {};
\end{tikzpicture} 
\end{minipage}
\vspace*{-0.2cm}
\caption{Sketch of geometrical setup (left) for test case in section \ref{sec:tc1:rot}, $\Omega_1$ is blue and $\Omega_2$ is green. And sketch of the time-independent part of the solution $U^1 / U^2$ (see \eqref{eq:tc1U12}) (right).}
\vspace*{-0.2cm}
\label{fig:tc1sketch}
\end{figure}
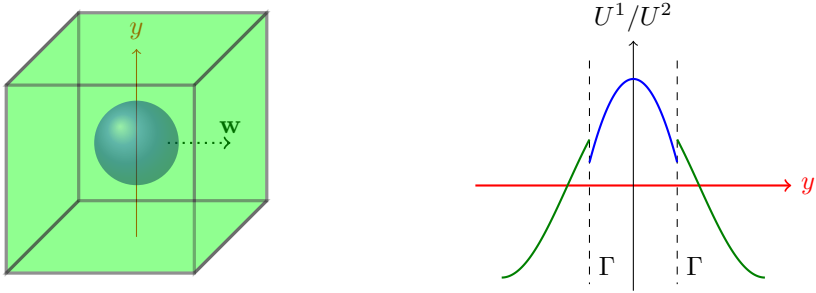
%
One phase is contained in the domain $\Omega_1(t)=\{\bx \in \Omega: \Vert x-(\bp_0 + r(t)\cdot \be_1) \Vert \leq R\}$, 
where $\bp_0$ is the center of the initial sphere and $r(t)$ the motion of the interface in $x_1$-direction, 
$\be_1$ is the corresponding unit vector. $R = \frac13$ is the radius of the sphere. 
The complementary domain is $\Omega_2(t)=\Omega \setminus \Omega_1(t)$. \\
The velocity field $\bw$ is given as $\bw = (\frac{\partial r}{\partial t}(t), 0, 0)$. As boundary conditions we apply suitable Dirichlet boundary conditions everywhere. \\
We prescribe these boundary conditions and the r.h.s. source term $f$, such that the solution is given by
$$
u(\bx,t) = \sin (k \pi t) \cdot U^{m}( \Vert \bx-(\bp_0+r(t)\cdot \be_1)\Vert ), \qquad \bx \in \Omega^{m}(t), \quad m=1,2
$$
with
\begin{equation} \label{eq:tc1U12}
U^1(y) = a+by^2
\quad \mbox{ and } \quad
U^2(y) = \cos(\pi y),
\end{equation}
where $a$ and $b$ are chosen s.t. the interface conditions hold. 
The diffusivities are $(\alpha_1,\alpha_2)=(10,20)$ and the Henry weights $(\beta_1,\beta_2)=(2,1)$ resulting in $a \approx 1.1569$ and $b \approx -8.1621$. The problem is considered in the time interval $[0,T]$ with $T=0.5$. We choose $\bp_0 = ( 0.5, 1 ,1)^T$ and $r(t)=\frac{1}{4 \pi} \sin(2 \pi t)$. For the approximation of the space-time interface we consider $m_s=m_t=1$.
We observe an error behaviour which is of (at least) second order in time and space ($\mathcal{O}(h^2+\Delta t^2)$) (see Figure \ref{fig:tc1_conv}). For the finest spatial resolution ($n_s=64$) we observe an order around $2.5$ for the convergence in time. In contrast to the previous test cases the spatial error dominates the overall error already for coarse temporal resolutions. We expect that for finer spatial resolutions and better geometry approximations ($m_t>1, m_s\geq 1$) one could retain the third order convergence in time. 
%

It is also relevant to study the accuracy of the method w.r.t. the interface condition. Therefore, in Table \ref{tbl:tc1:rot:linmove} we consider the error $$\Vert \nu^{\frac12} [\beta u_h] \Vert_{L^2(\Gamma_{\ast,h})} = \left( \int_0^T \int_{\Gamma_h(t)} [\beta u_h]^2 \, ds \, dt \right)^{\frac12}$$ under space and time refinement and also observe an $\mathcal{O}(h^2+\Delta t^2)$ behaviour. 
\begin{table}[ht]
\small
\begin{center}
  \begin{tabular}{|r|r|r|r||r|r|}
\hline
$\nt$ $\backslash$ $\ns$ & 8 & 16 & 32 & \ccg 64 & \ccg \eoct \\
\hline
2 & 0.0495 & 0.00700 & 0.0198 & \ccg   0.0587 & \ccg \\
4 & 0.0430 & 0.00384 & 0.00567 & \ccg   0.0227 & \ccg 1.37 \\
8 & 0.0417 & 0.00253 & 0.00164 & \ccg  0.00517 & \ccg 2.13 \\
16 & 0.0414 & 0.00205 & 0.000716 & \ccg  0.00117 & \ccg 2.14 \\
32 & 0.0413 & 0.00190 & 0.000523 & \ccg 0.000275 & \ccg 2.10 \\
\hline
\ccg 64 & \ccg 0.0413 & \ccg 0.00186 & \ccg 0.000477 & \ccg 0.000131 & \ccg 1.07 \\
\hline
\ccg \eocs & \ccg & \ccg 4.48 & \ccg 1.96 & \ccg 1.86 & \ccg \\
\hline
\end{tabular}
\end{center}
\caption{Interface error $\Vert \nu^{\frac12} [\beta u_h] \Vert_{L^2(\Gamma_{\ast,h})}$ for different refinements in time and space for the test case in section \ref{sec:tc1:rot}. The last column shows the estimated order of convergence w.r.t. time (\eoct) using the finest spatial resolution ($\ns=64$), the last row shows the estimated order of convergence w.r.t. space (\eocs) using the finest temporal resolution ($\nt=64$).}
\label{tbl:tc1:rot:linmove}
\vspace*{-0.5cm}
\end{table}
%

%

\subsection{Deforming bubble in a vortex}
\label{sec:tc7:vortex}
\begin{figure}[ht]
  \begin{minipage}{0.166\textwidth}
    \includegraphics[width=1.2\textwidth]{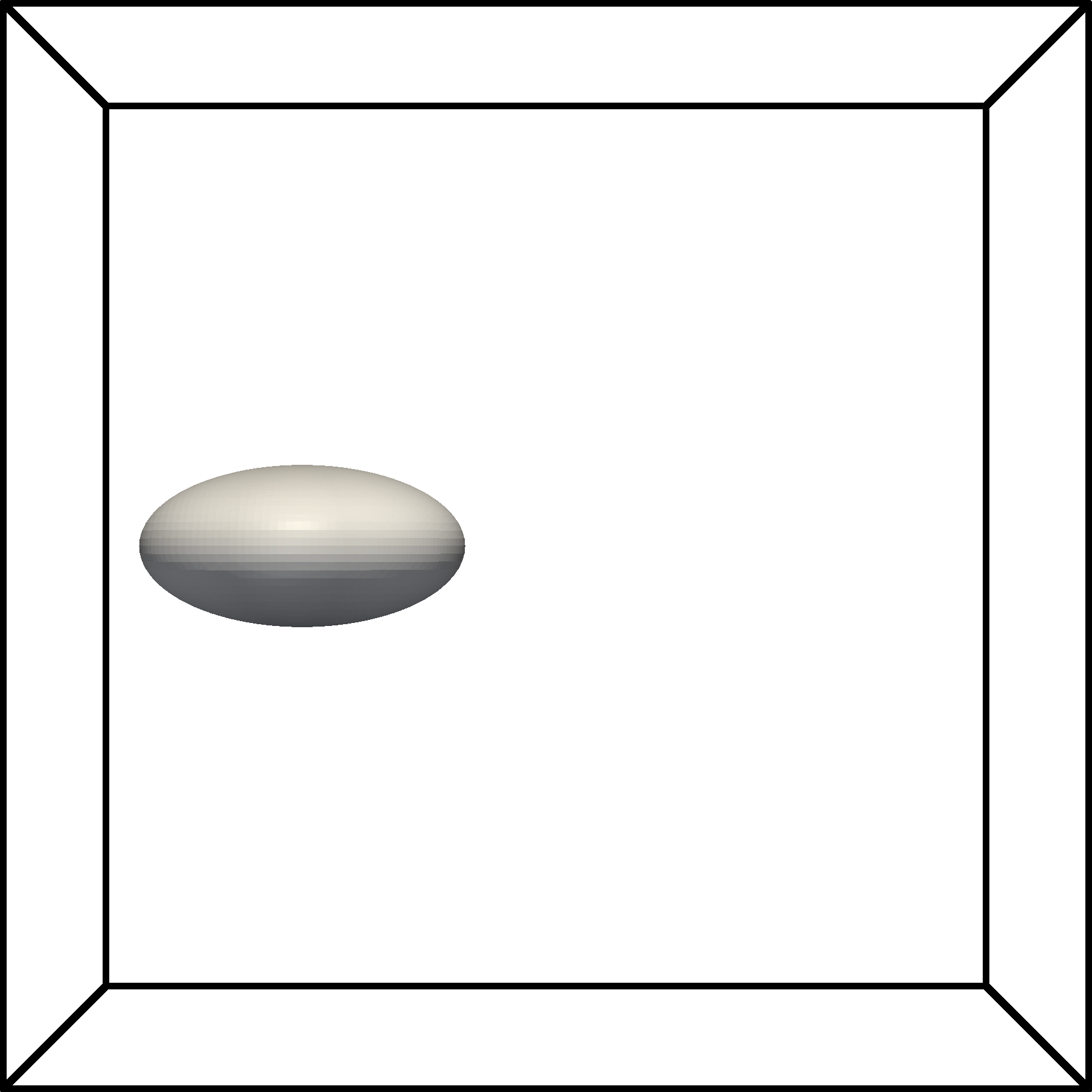} 
    \centering{$\,$\hspace*{0.1cm} \small $\quad t= 0$}
  \end{minipage}
  $\ \ $
  \begin{minipage}{0.166\textwidth}
    \includegraphics[width=1.2\textwidth]{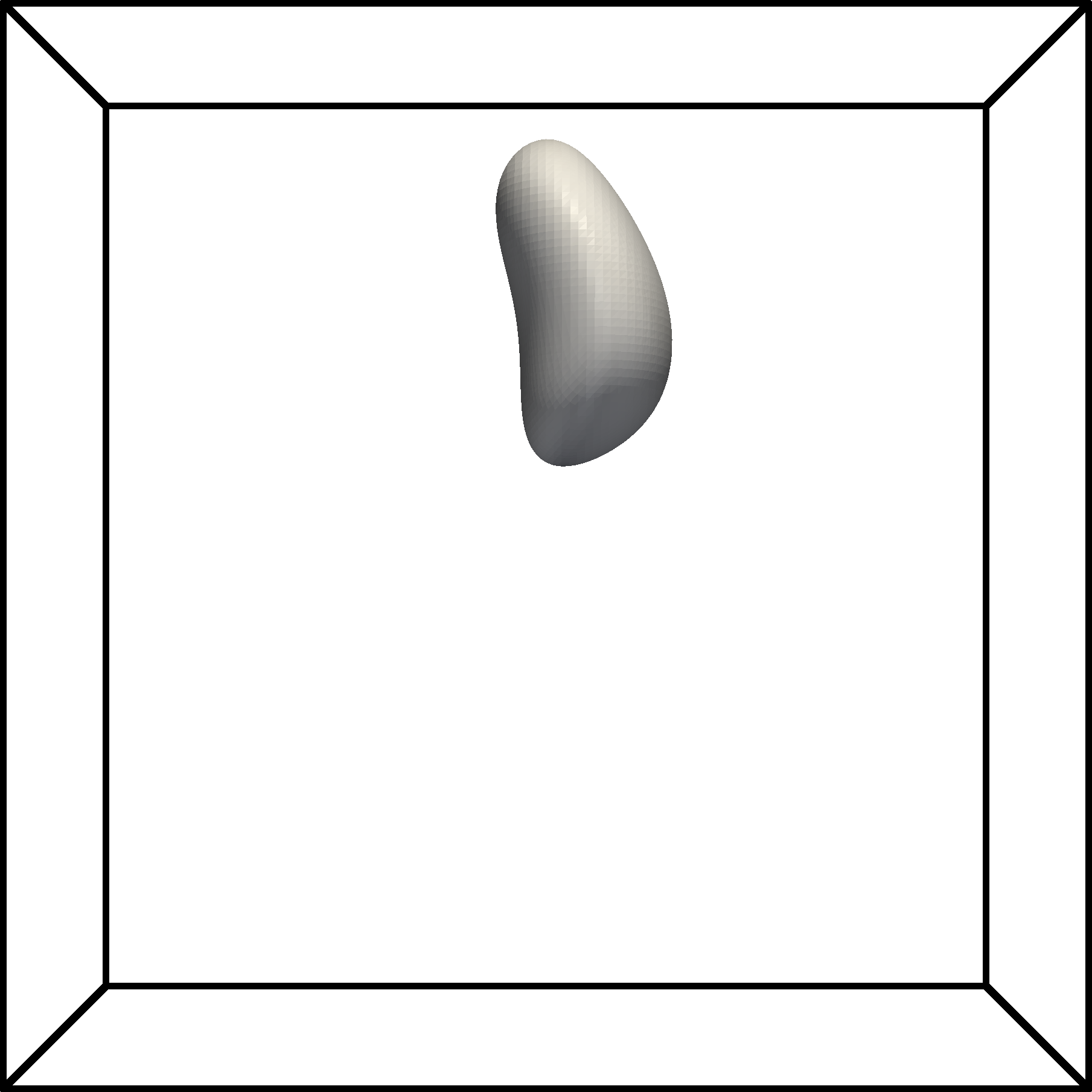}
    \centering{$\,$\hspace*{0.1cm} \small $\quad t= 5$}
  \end{minipage}
  $\ \ $
  \begin{minipage}{0.166\textwidth}
    \includegraphics[width=1.2\textwidth]{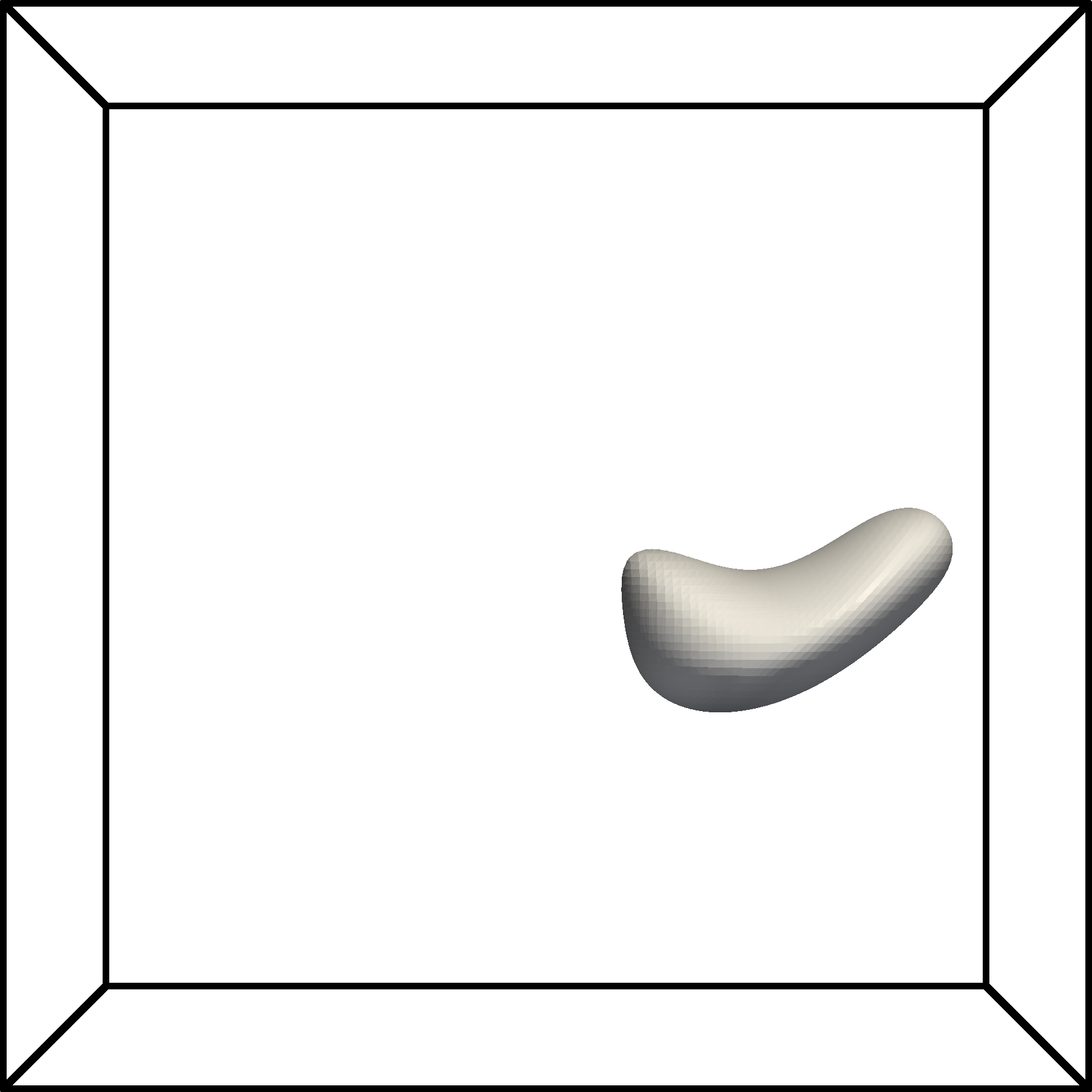}
    \centering{$\,$\hspace*{0.1cm} \small $\quad t= 10$}
  \end{minipage}
  $\ \ $
  \begin{minipage}{0.166\textwidth}
    \includegraphics[width=1.2\textwidth]{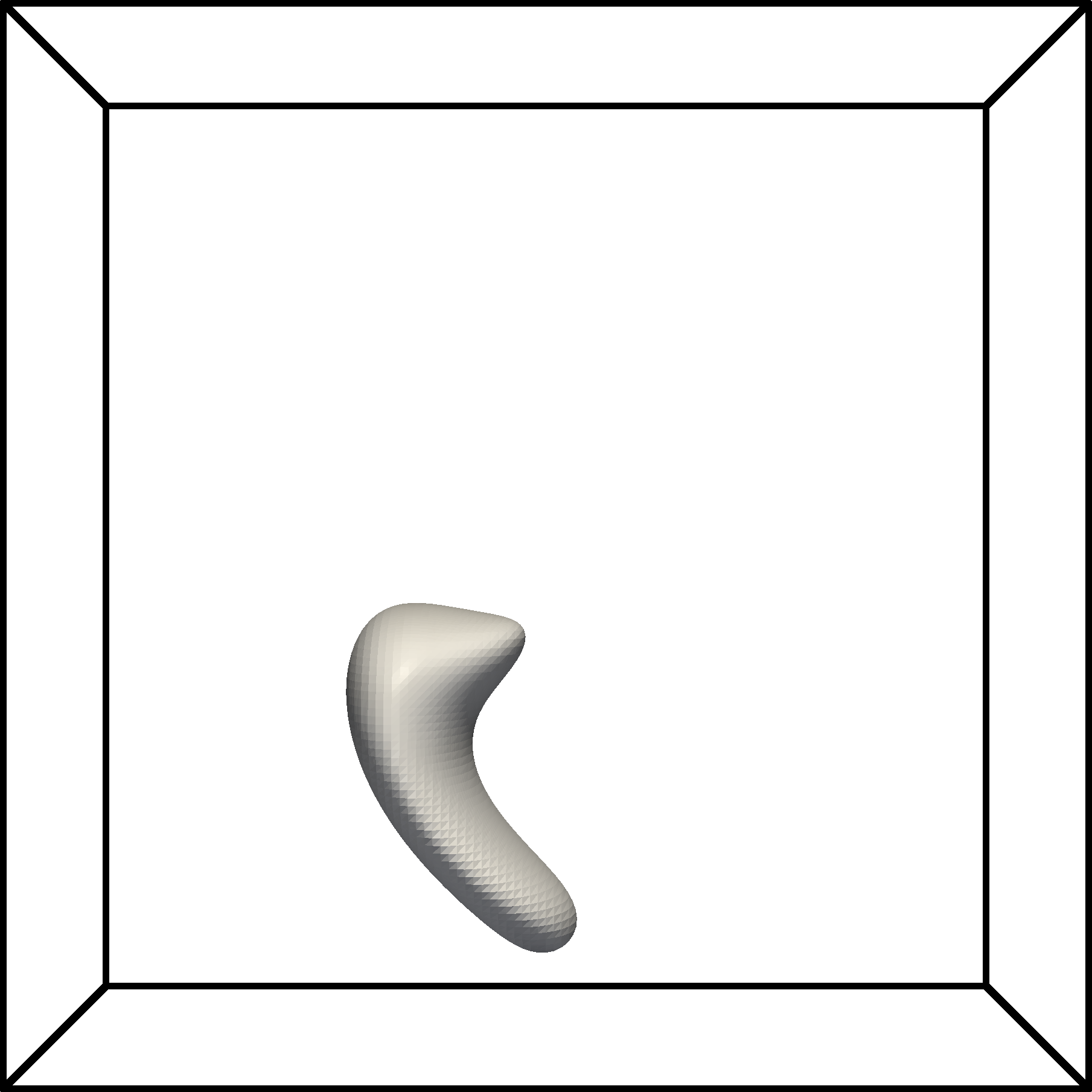}
    \centering{$\,$\hspace*{0.1cm} \small $\quad t= 15$}
  \end{minipage}
  $\ \ $
  \begin{minipage}{0.166\textwidth}
    \includegraphics[width=1.2\textwidth]{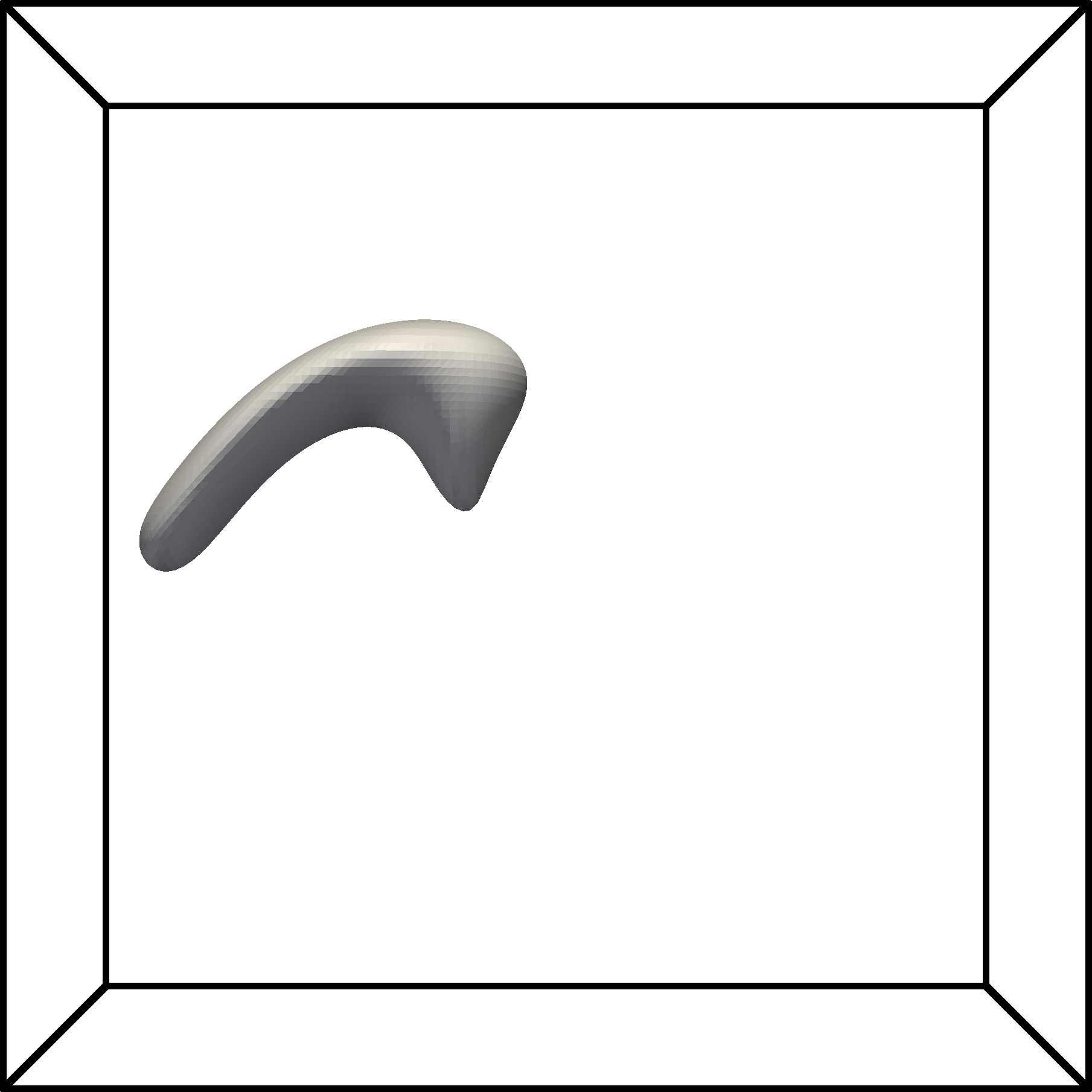}
    \centering{$\,$\hspace*{0.1cm} \small $\quad t= 20$}
  \end{minipage}
  $\ \ $
  \vspace*{-0.2cm}
  \caption{Interface position at several times for the test case in section \ref{sec:tc7:vortex}, view from the top.}\label{fig:tc7_vort_if}
\end{figure}

As a last example we consider a more complex configuration. We consider an ellipsoidal bubble which is deforming under a vortex velocity field. The domain is $\Omega = [0,2] \times [0,2] \times [0,1]$ and the velocity field is given as:
\begin{equation}
 \bw(\bx,t) = \bw(\bx) = q(r(\bx)) 
\cdot ( x_2 - c_2, c_1 - x_1, 0 )^T
\end{equation} 
where $(c_1,c_2,x_3)$ with $c_1=c_2=1$ describes the rotation axes of the vortex and
$$
 r(\bx) := \sqrt{(x_1-c_1)^2 + (x_2-c_2)^2}.
$$ 
The angular velocity varies with a changing distance to the rotation axes. The dependency is described by $q(r)$ which is chosen as follows:
\begin{equation}
q(r) = \frac{\pi}{10} \cdot \left\{ 
\begin{array}{cr@{}c@{}l}
1 + 2^2 \frac{5^{10}}{3^{13}} \, r^2 (0.9-r)^{3} &  & r & \leq 0.9 \\
10^4 (x-0.8)^2 (x-1)^2 & 0.9 < & r & < 1 \\
0 & 1 \leq & r &
  \end{array}
\right.
\end{equation}
The time interval is $[0,T]$ with $T=20$ which corresponds to one full rotation of the bubble. 
Due to $q(r) \neq \mathrm{const}$ the bubble deforms during that process. In Figure \ref{fig:tc7_vort_if} the interface at different time levels is shown. 
%

Initially we set 
$$
\Omega_1(0) = \left\{ 3^2 (x_1-b_1^0)^2 + 6^2 (x_2-b_2^0)^2 + 3^2 (x_3-b_3^0)^2 \leq 1 \right\}
$$
with $(b_1^0,b_2^0,b_3^0) = (0.5, 1 , 0.5)$, the center of the initial bubble. The complementary initial domain is $\Omega_{2}(0) = \Omega \setminus \Omega_{1}(0)$. The interface evolves along the characteristics given by $\dot{\bx}(t) = \bw(\bx,t)$ which can be calculated explicitely. The diffusivities are $(\alpha_1, \alpha_2) = (8 \cdot 10^{-3}, 4 \cdot 10^{-3})$ and the Henry weights are $(\beta_1,\beta_2) = (1, 1.75)$. As initial condition we choose $u_1(\bx,0) = 1$ and $u_2(\bx,0) = 0$. Note that $u(\cdot,0)$ does not fulfill the interface condition \eqref{eq3} which leads to a parabolic boundary layer of size $\mathcal{O}(\sqrt{\alpha t})$. 
\begin{figure}[ht]
  \begin{minipage}{0.456\textwidth}
    \includegraphics[width=0.9\textwidth]{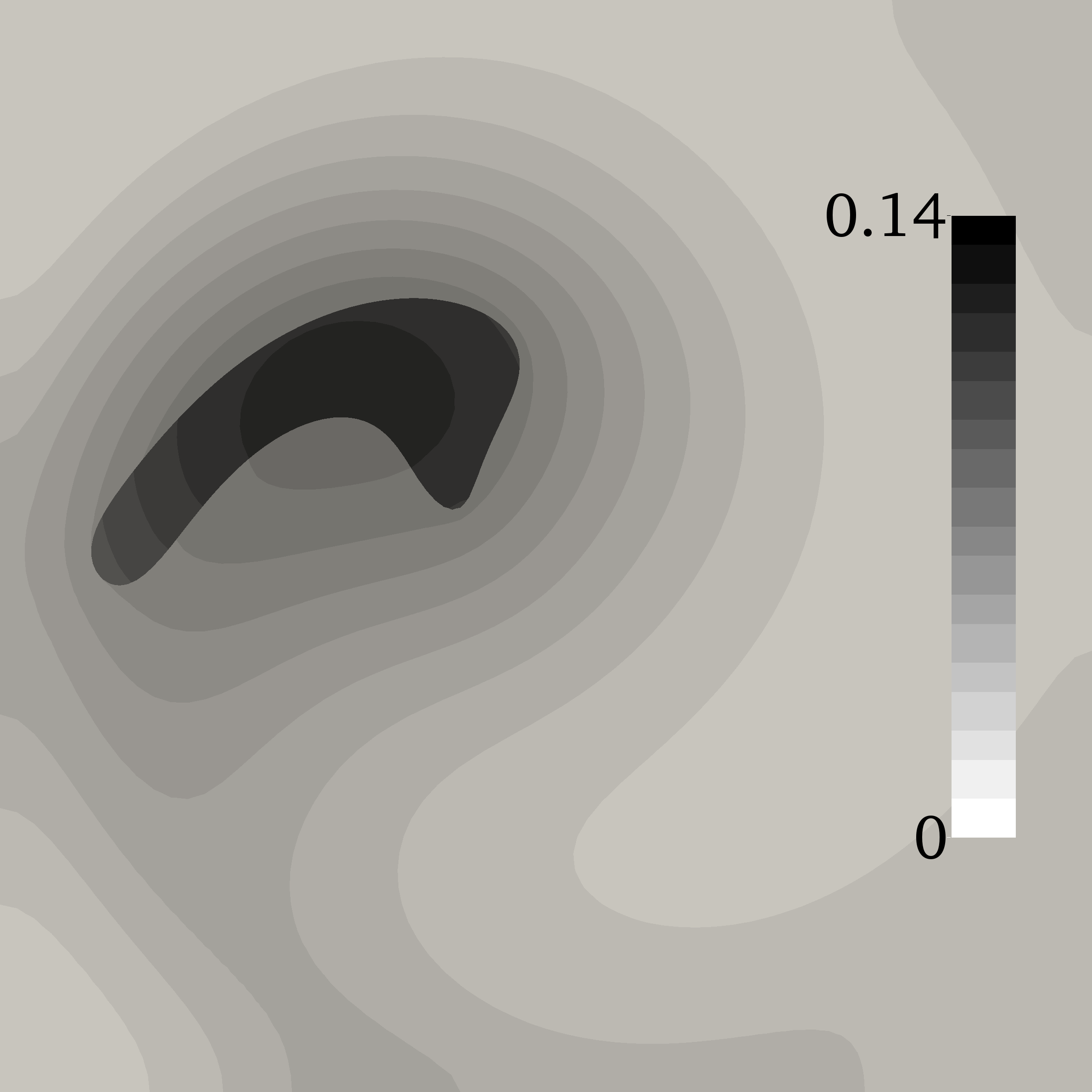} 
  \end{minipage}
  \begin{minipage}{0.5\textwidth}
    \begin{center}
      \begin{tabular}{|r|r|r|}
        \hline
        $\nt$ & $\Vert u-u_{\mathrm{ref}}(\cdot, T) \Vert_{L^2(\Omega)} $ & \eoct \\
        \hline
         32 & $1.22569 \cdot 10^{-3}$ &  \\
         64 & $3.32611 \cdot 10^{-4}$ & 1.88 \\
        128 & $9.65349 \cdot 10^{-5}$ & 1.78 \\
        256 & $4.08119 \cdot 10^{-5}$ & 1.24 \\
        \hline
      \end{tabular}
    \end{center}
  \end{minipage}
  \vspace*{-0.2cm}
  \caption{Numerical solution and convergence table of the test case in section \ref{sec:tc7:vortex}. The contour plot (left) shows the numerical solution ($n_s = 96$, $n_t= 1024$) of the test case in section \ref{sec:tc7:vortex} at $T=20$ in the cutting plane $z=0.5$. The table shows the convergence in time w.r.t. a reference solution.}\label{fig:tc7_vort_sol}
\end{figure}
%

For this setup we do not know the exact solution.  We therefore computed a reference solution $u_{\mathrm{ref}}$ on a spatial fine mesh ($96 \times 96 \times 48$) with a characteristic mesh size of $h= 1/96$ and a temporal resolution of $1024$ timesteps for the whole time interval, i.e. $\Delta t_{\mathrm{ref}} \approx 2 \cdot 10^{-2}$. This solution is denoted by $u_{\mathrm{ref}}$. We only investigate the temporal convergence by varying the number of time steps $n_t$. For the approximation of the space-time interface we consider $m_s=2$, $m_t=1$. A contour plot of the solution and the convergence table is given in Figure \ref{fig:tc7_vort_sol}. We observe a convergence order in time of nearly two. 
The range in which we observe this order of convergence is limited by $n_t < 256$. This is due to several effects. At first the reference solution is different from the exact solution. Secondly due to the dependence of the finite element space and the approximated space-time interface on the time step size the impact of the spatial discretization errors can be different for different time step sizes. Hence $u-u_{\mathrm{ref}}$ does contain spatial errors which do not vanish for $\Delta t \rightarrow 0$.


%
\section{A strategy to decompose intersected 4-prisms into pentatopes}\label{sec:decomp4d}
In this section we introduce a decomposition strategy that allows for a decomposition of four dimensional prisms into pentatopes as needed in section \ref{sec:calcints}. This approach is new.

Firstly, we introduce the definitions of relevant four dimensional geometries in section \ref{sec:prims4d}. 
The decomposition of a 4-prism into four pentatopes is presented in section \ref{sec:decomprefprism}. This is already needed to construct (via interpolation of the level-set function) the piecewise planar space-time interface in section \ref{sec:stifapprox}. 
\begin{figure}[ht]
  \begin{center}
    \tikzstyle{decision} = [diamond, draw, fill=blue!20, 
    text width=4.5em, text badly centered, inner sep=0pt]

    \tikzstyle{cutblock} = [rectangle, draw, fill=blue!10, 
    text width=5em, text centered, rounded corners, minimum height=4em]

    \tikzstyle{mixedblock} = [rectangle, draw, fill=blue!10, 
    text width=5em, text centered, rounded corners, minimum height=4em]

    \tikzstyle{uncutblock} = [rectangle, draw, fill=blue!10, 
    text width=5em, text centered, rounded corners, minimum height=4em]

    \tikzstyle{line} = [draw, -latex', thick]
    \tikzstyle{cloud} = [draw, ellipse,fill=red!20, node distance=3cm,
    minimum height=2em]
    \begin{tikzpicture}[node distance = 1cm, auto]
      \node [cutblock] (init) {one cut 4-prism };
      \node [mixedblock, below =1cm of init] (afterinit) {four pentatopes};
      \node [cutblock, below = 1cm of afterinit] (cutpenta) {cut pentatope \\ \small (sec. \ref{sec:decompintersectpenta})};
      \node [uncutblock, above right =1cm and 8cm of afterinit] (uncutpenta) {uncut pentatope};
      \node [cutblock, above right = 0.25cm and 1cm of cutpenta] (case1) {case 1 \\ \small(sec. \ref{sec:case1})};
      \node [cutblock, below right =-0.25cm and 1cm of cutpenta] (case2) {case 2 \\ \small(sec. \ref{sec:case2})};
      \node [uncutblock, above right=-0.5cm and 5cm of case1] (case1a) {one uncut pentatope};
      \node [uncutblock, below right=   0cm and 1cm of case1] (case1b) {one 4-prism};
      \node [uncutblock, right= 2cm of case1b] (case1btopenta) {four uncut pentatopes};
      \node [uncutblock, below right=-0.5cm and 1cm of case2] (case2b) {one hypertriangle};
      \node [uncutblock, right= 2cm of case2b] (case2btopenta) {six uncut pentatopes};
      \path [line] (init) -- node {\small(sec. \ref{sec:decomprefprism})} (afterinit);
      \path [line] (afterinit) -- (cutpenta);
      \path [line] (afterinit) -- (uncutpenta);
      \path [line] (cutpenta) -- (case1);
      \path [line] (cutpenta) -- (case2);
      \path [line] (case1) -- (case1a);
      \path [line] (case1) -- (case1b);
      \path [line] (case2) -- (case1b);
      \path [line] (case1b) -- node {\small(sec. \ref{sec:decomprefprism})} (case1btopenta);
      \path [line] (case2) -- (case2b);
      \path [line] (case2b) -- node {\small(sec. \ref{sec:decomprefhypertrig})} (case2btopenta);


      \begin{pgfonlayer}{background}
        \path (uncutpenta.east |- uncutpenta.north)+(+0.32,+0.32) node (a) {};
        \path (case2btopenta.south -| case2btopenta.west)+(-0.32,-0.82) node (b) {};
        \path (case2b.west |- case2b.south)+(-0.32,-0.82) node (c) {};
        \path (case1b.north -| case1b.east)+(+0.32,+0.32) node (d) {};
        \path (init.west |- init.north)+(-0.32,0.32) node (e) {};
        \path (case2b.west |- case2b.south)+(-0.72,-0.82) node (f) {};

        \path (f)+(-0.25,0.25) node (f1) {};
        \path (f)+(-5.35,1) node (f2) {};

        \path (c)+(0.25,0.25) node (c1) {};
        \path (c)+(2.35,1) node (c2) {};

        \path (b)+(0.25,0.25) node (b1) {};
        \path (b)+(2.35,1) node (b2) {};

        \path[fill=green!20,rounded corners, draw=black!50, dashed]
        (a) rectangle (b);           
        \path[fill=yellow!20,rounded corners, draw=black!50, dashed]
        (c) rectangle (d);           
        \path[fill=orange!20,rounded corners, draw=black!50, dashed]
        (e) rectangle (f);           
      \end{pgfonlayer}

      \path (f2) +(2.71,-0.6) node (aasdf) {\small cut geometries};
      \path (c2) +(-1.04,-0.6) node (aasdf2) {\small uncut geometries};
      \path (b2) +(-1.041,-0.6) node (aasdf3) {\small uncut pentatopes};



      
    \end{tikzpicture}
  \end{center}
  \vspace*{-0.2cm}
  \caption{Algorithmic structure of the decomposition strategy proposed in section \ref{sec:decomp4d}.}
  \label{fig:algstruct}
  \vspace*{-0.1cm}
\end{figure}
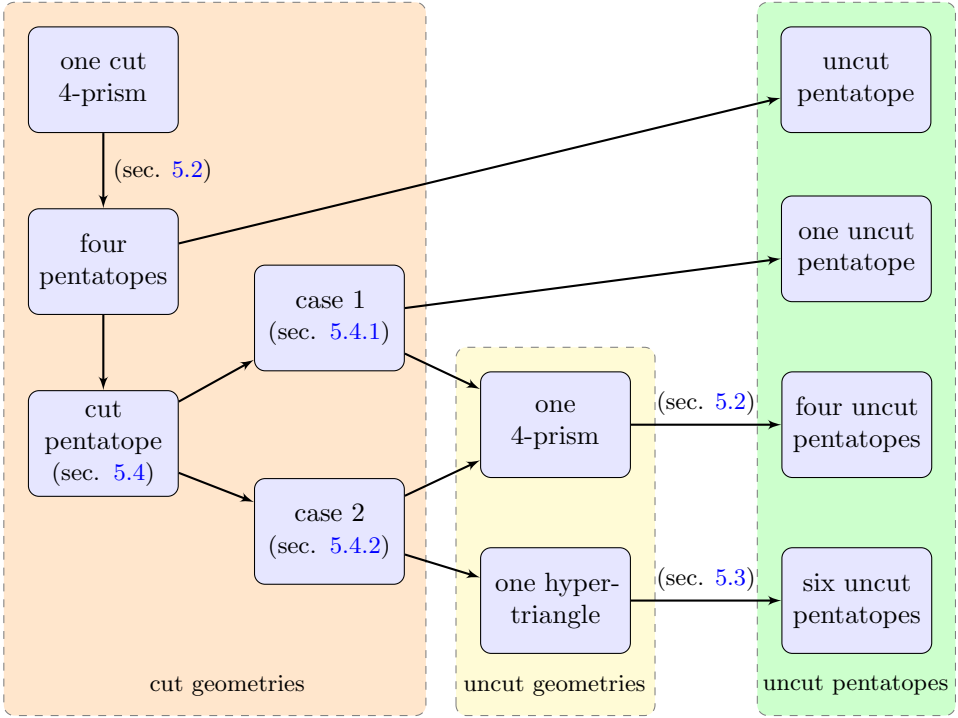
%
In section \ref{sec:decompintersectpenta} a strategy is presented that allows us to decompose a pentatope which is intersected by a hyperplane (representing an approximation of the space-time interface) into pentatopes which are not intersected. Figure \ref{fig:algstruct} sketches the algorithmic structure of the decomposition strategy.
In this algorithm we need a particular geometrical object, that we call \emph{hypertriangle}, which can be decomposed into six pentatopes following the decomposition rule in section \ref{sec:decomprefhypertrig}.
\begin{remark}[Small angles]
The resulting pentatopes/tetrahedra in this decomposition can have arbitrary small angles. Note that this does not lead to stability problems as we are using the decomposition only for numerical integration.
\end{remark}
\subsection{Definition of simple geometries in four dimensions}
\label{sec:prims4d}
By $\be^i \in \rr^n$ we denote the $i$-th unit vector with $(\be^i)_j = \delta_{i,j}$ for $i=1,..,n$ and $\be^0:=0$. 
\begin{definition}[4-simplex / pentatope]
  Let $\bx^i \in \rr^4$ for $i =1,..,5$ and $\mathbf{d}^{i,j} := \bx^i - \bx^j$. Iff the vectors $\mathbf{d}^{i,1}$ for $ i = 2,..,5 $ are linearly independent, we call the convex hull $\mP = \textbf{conv}( \{\bx^i\}_{i=1,..,5})$ the \emph{4-simplex} or \emph{pentatope}.
\end{definition}
\begin{remark}[reference pentatope]
\rm
Every pentatope $\mP$ can be represented as an affine transformation applied to the reference pentatope
$\hat{\mP} = \textbf{conv}(\{\be^i\}_{i=0,..4})$. The transformation has the form 
\vspace*{-0.05cm}
$$
\Phi: \hat{\mP} \rightarrow \mP, (\hat{x}_1,\hat{x}_2,\hat{x}_3,\hat{x}_4) \rightarrow 
 \sum_{i=1}^5 \hat{\lambda}_i \bx^i,
\vspace*{-0.05cm}
$$
where $\hat{\lambda}_i(\hat{x}_1,\hat{x}_2,\hat{x}_3,\hat{x}_4)$ are the barycentric coordinates of $\hat{\mP}$ w.r.t. vertex $\be^{i-1}$.
\end{remark}
\begin{definition}[4-prism]
  Let $\bx^i \in \mathbb{R}^4$ for $i = 1,..,4$ and $\mathbf{y} \in \mathbb{R}^4$. Iff $\{\bx^i\}_{i=1,..,4}$ defines a 3-simplex (tetrahedron) $T=\textbf{conv}(\{\bx^i\}_{i=1,..,4})$ and $\mathbf{y}$ is linearly independent of $\{ \mathbf{d}^{i,1} \}_{i=2,..,4}$, with $\mathbf{d}^{i,j} := \bx^i - \bx^j$, the set
$$ \mQ = \bconv(\{ \bx^i \}_{i=1,..,4},\{ \bx^i + \mathbf{y} \}_{i=1,..,4}) =  \{ \bx + \alpha \mathbf{y}, \bx \in \bconv(\{\bx^i\}_{i=1,..,4}), \alpha \in [0,1] \}$$ 
is called \emph{4-prism}.
\end{definition}
\begin{remark}[reference 4-prism]
\rm
Every 4-prism can be represented as an affine linear transformation applied to the reference 4-prism
$\hat{\mathcal{Q}} = \textbf{conv}(\{\be^i\}_{i=0,..,3} \},\{\be^i + \be^4\}_{i=0,..,3} \})$.
The transformation has the form 
$$
\Phi: \hat{\mQ} \rightarrow \mQ, (\hat{x}_1,\hat{x}_2,\hat{x}_3,\hat{x}_4) \rightarrow 
\sum_{i=1}^4 \hat{\mu}_i \bx^i + \hat{x}_4 \by,
$$
where $\hat{\mu}_i(\hat{x}_1,\hat{x}_2,\hat{x}_3)$ are the barycentric coordinates of the reference tetrahedron $\hat{T} = \bconv(\{\be^i\}_{i=0,..3})$.
\end{remark}
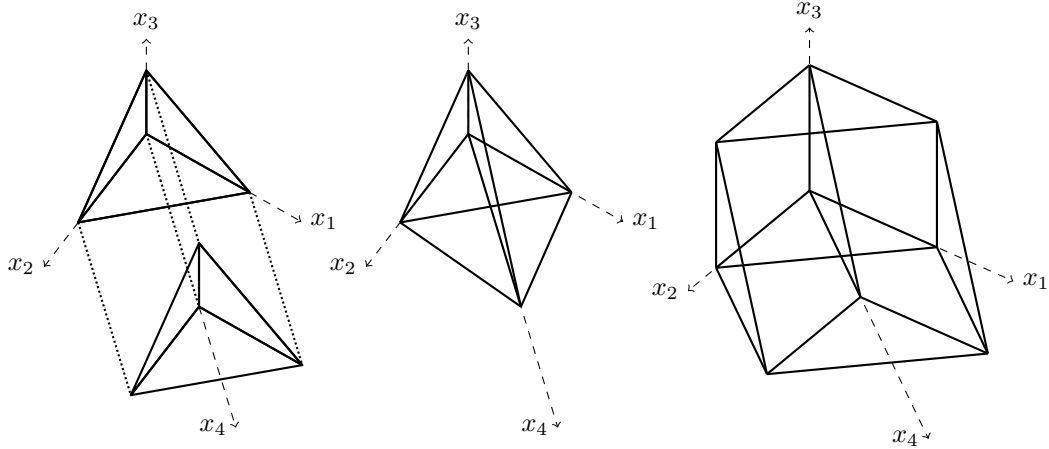
\begin{figure}[ht]
  \vspace*{-0.3cm}
  \hspace*{-0.45cm}
  \begin{minipage}{0.33\textwidth}  
    \begin{tikzpicture}[line join=round]
      \tikzstyle{transparent cone} = [fill=blue!20,fill opacity=0.8]\filldraw[thick,fill=lightgray,fill opacity=0.0](0,0)--(1.373,-.777)--(-.906,-1.177)--cycle;
      \draw[thick,densely dotted,draw=black](0,0)--(.187,-.613);
      \filldraw[thick,fill=lightgray,fill opacity=0.0](0,0)--(-.906,-1.177)--(0,.845)--cycle;
      \filldraw[thick,fill=lightgray,fill opacity=0.0](0,0)--(0,.845)--(1.373,-.777)--cycle;
      \draw[dashed,draw=black,arrows=-](0,0)--(1.373,-.777);
      \draw[dashed,draw=black,arrows=-](0,0)--(0,.845);
      \draw[dashed,draw=black,arrows=-](0,0)--(-.906,-1.177);
      \filldraw[thick,fill=lightgray,fill opacity=0.0](-.906,-1.177)--(1.373,-.777)--(0,.845)--cycle;
      \draw[dashed,draw=black,arrows=->](1.373,-.777)--(2.059,-1.166);
      \draw[thick,densely dotted,draw=black](1.373,-.777)--(2.073,-3.075);
      \draw[dashed,draw=black,arrows=->](-.906,-1.177)--(-1.359,-1.766);
      \draw[thick,densely dotted,draw=black](-.906,-1.177)--(-.206,-3.475);
      \draw[thick,densely dotted,draw=black](.187,-.613)--(.7,-2.298);
      \draw[thick,densely dotted,draw=black](0,.845)--(.7,-1.453);
      \draw[dashed,draw=black,arrows=->](0,.845)--(0,1.268);
      \filldraw[thick,fill=lightgray,fill opacity=0.0](.7,-2.298)--(2.073,-3.075)--(-.206,-3.475)--cycle;
      \draw[dashed,draw=black,arrows=->](.7,-2.298)--(1.19,-3.906);
      \filldraw[thick,fill=lightgray,fill opacity=0.0](.7,-2.298)--(-.206,-3.475)--(.7,-1.453)--cycle;
      \filldraw[thick,fill=lightgray,fill opacity=0.0](.7,-2.298)--(.7,-1.453)--(2.073,-3.075)--cycle;
      \filldraw[thick,fill=lightgray,fill opacity=0.0](.7,-1.453)--(2.073,-3.075)--(.7,-1.453)--cycle;
      \path (2.059,-1.166) node[right] {$x_1$}
      (0,1.268) node[above] {$x_3$}
      (-1.359,-1.766) node[left] {$x_2$}
      (1.19,-3.906) node[left] {$x_4$};\end{tikzpicture}
  \end{minipage}
  \hspace*{-0.25cm}
  \begin{minipage}{0.33\textwidth}  
    \begin{tikzpicture}[line join=round]
      \tikzstyle{transparent cone} = [fill=blue!20,fill opacity=0.8]\draw[thick,draw=black](0,0)--(.7,-2.298);
      \draw[thick,draw=black](0,0)--(0,.845);
      \draw[thick,draw=black](0,0)--(-.906,-1.177);
      \draw[thick,draw=black](0,0)--(1.373,-.777);
      \draw[dashed,draw=black,arrows=->](0,0)--(2.059,-1.166);
      \draw[dashed,draw=black,arrows=->](0,0)--(0,1.268);
      \draw[dashed,draw=black,arrows=->](0,0)--(-1.359,-1.766);
      \draw[dashed,draw=black,arrows=->](0,0)--(1.19,-3.906);
      \draw[thick,draw=black](1.373,-.777)--(.7,-2.298);
      \draw[thick,draw=black](1.373,-.777)--(-.906,-1.177);
      \draw[thick,draw=black](1.373,-.777)--(0,.845);
      \draw[thick,draw=black](-.906,-1.177)--(0,.845);
      \draw[thick,draw=black](-.906,-1.177)--(.7,-2.298);
      \draw[thick,draw=black](0,.845)--(.7,-2.298);
      \path (2.059,-1.166) node[right] {$x_1$}
      (0,1.268) node[above] {$x_3$}
      (-1.359,-1.766) node[left] {$x_2$}
      (1.19,-3.906) node[left] {$x_4$};\end{tikzpicture}
  \end{minipage}
  \hspace*{-0.25cm}
  \begin{minipage}{0.33\textwidth}  
    \begin{tikzpicture}[line join=round]
      \tikzstyle{transparent cone} = [fill=blue!20,fill opacity=0.8]\draw[thick,draw=black](0,0)--(.677,-1.415);
      \draw[thick,draw=black](0,0)--(0,1.671);
      \draw[thick,draw=black](0,0)--(-1.242,-1.026);
      \draw[thick,draw=black](0,0)--(1.694,-.753);
      \draw[dashed,draw=black,arrows=->](0,0)--(2.71,-1.204);
      \draw[dashed,draw=black,arrows=->](0,0)--(0,2.172);
      \draw[dashed,draw=black,arrows=->](0,0)--(1.581,-3.301);
      \draw[dashed,draw=black,arrows=->](0,0)--(-1.615,-1.334);
      \draw[thick,draw=black](1.694,-.753)--(2.371,-2.167);
      \draw[thick,draw=black](1.694,-.753)--(1.694,.918);
      \draw[thick,draw=black](1.694,-.753)--(-1.242,-1.026);
      \draw[thick,draw=black](0,1.671)--(.677,-1.415);
      \draw[thick,draw=black](0,1.671)--(-1.242,.645);
      \draw[thick,draw=black](0,1.671)--(1.694,.918);
      \draw[thick,draw=black](-1.242,-1.026)--(-1.242,.645);
      \draw[thick,draw=black](-1.242,-1.026)--(-.565,-2.441);
      \draw[thick,draw=black](1.694,.918)--(2.371,-2.167);
      \draw[thick,draw=black](1.694,.918)--(-1.242,.645);
      \draw[thick,draw=black](-1.242,.645)--(-.565,-2.441);
      \draw[thick,draw=black](.677,-1.415)--(-.565,-2.441);
      \draw[thick,draw=black](.677,-1.415)--(2.371,-2.167);
      \draw[thick,draw=black](2.371,-2.167)--(-.565,-2.441);
      \path (2.71,-1.204) node[right] {$x_1$}
      (0,2.172) node[above] {$x_3$}
      (-1.615,-1.334) node[left] {$x_2$}
      (1.581,-3.301) node[left] {$x_4$};\end{tikzpicture}
  \end{minipage}
  \vspace*{-0.3cm}
  \caption{Sketch of reference geometries. Reference 4-prism $\hat{\mathcal{Q}}$ (left), reference pentatope $\hat{\mP}$ (center) and reference hypertriangle $\hat{\mathcal{H}}$ (right). The dotted line in the left picture are parallel to the $x_4$-axes and connect the tetrahedra at $x_4=0$ and $x_4=1$. }\label{fig:refsketches}
\end{figure}
%
The next geometry is a little bit more complex. It later occurs as one part of a pentatope cut by a hyperplane. 

\begin{definition}[hypertriangle]\label{def:hypertriangle}
We define the reference hypertriangle as 
\begin{eqnarray*}
\hat{\mH} &:=& \{(x_1,x_2,x_3,x_4)\in \rr_+^4, x_1+x_2\leq 1, x_3 +x_4 \leq 1\} \\
 &=& \bconv(\{ \hat{\bx}^{i,j} \}_{i=1,..,3,j=1,..,3}) = \hat{\mathcal{K}} \times \hat{\mathcal{K}}
\end{eqnarray*}
where $\hat{\mathcal{K}} \subset \rr^2$ denotes the reference triangle $ \hat{\mathcal{K}} = \bconv(\{\chi^1,\chi^2,\chi^3\}) \subset \rr^2 $ with $\chi^i=\be^{i-1} \in \rr^2$, $i=1,2,3$ and $\hat{\bx}^{i,j}=(\chi^i,\chi^j) \in \rr^4$, $i,j=1,2,3$.
Now, let $\bx^{i,j} \in \rr^4,\ i,j=1,2,3$. The convex hull $\mathcal{H}=\bconv(\{ \bx^{i,j}\}_{i,j=1,2,3})$ is called a \emph{hypertriangle} iff there exists a transformation 
$$
\Phi: \hat{\mathcal{H}} \rightarrow \rr^4, (\hat{x}_1,\hat{x}_2,\hat{x}_3,\hat{x}_4) \rightarrow \sum_{i=1}^3 \sum_{j=1}^3 \hat{\rho}_i(\hat{x}_1,\hat{x}_2) \hat{\rho}_j(\hat{x}_3,\hat{x}_4) \bx^{i,j},
$$
where $\hat{\rho}_i(\hat{x}_1,\hat{x}_2)$ is the barycentric coordinate of the reference triangle $\hat{\mathcal{K}}$corresponding to the vertex $\chi^i$ such that there holds $\Phi(\hat{\mH}) = \mH$.
\end{definition}

\begin{remark}[Sketches] \label{rem:sketches}
The sketches in Figures \ref{fig:refsketches}, \ref{fig:decomprefprism4} and \ref{fig:decomprefhypertrig} in this section show two dimensional parallel projections of four dimensional objects. Straight lines in the sketch represent a line (in four dimensions) between two vertices. Note that the preimage of a point in the two dimensional sketch of the parallel projection is a two-dimensional set. 
\end{remark}

\subsection{Decomposition of a 4-prism into four pentatopes} \label{sec:decomprefprism}

We consider an arbitrary prism element $\mathcal{Q}^T = T \times I_n $ with a tetrahedral element $T$ and a time interval $I_n$. 
For each $\mQ^T$ there exists a linear transformation $\Phi$ mapping from the reference 4-prism $\hat{\mQ}$ to $\mQ^T$ which is of the form $\Phi(\hat{\bx},\hat{t}) = ( \Phi_x(\hat{\bx}), \Phi_t(\hat{t}))^T$ with the time transformation $\Phi_t(\hat{t}) = \hat{t} \cdot t_n +(1-\hat{t}) \cdot t_{n-1}$ and the space transformation $\Phi_x(\hat{\bx})$ mapping from the reference tetrahedron $\hat{T}$ to $T$.

It is sufficient to consider the decomposition of the reference 4-prism $\hat{\mQ}$ into four pentatopes as applying $\Phi$ to each pentatope of this decomposition results in a valid decomposition of $\mathcal{Q}$ into four pentatopes.
With $\pbx{i} := \be^{i-1}$ and $\pby{i} := \be^{i-1} + \be^4$ for $i = 1,..,4$ for the reference 4-prism there holds 
$\hat{\mQ} = \bconv( \{\pbx{i}\}_{i=1,..,4}, \{\pby{i}\}_{i=1,..,4})$. 
We decompose $\hat{\mathcal{Q}}$ into four pentatopes $\hat{\mP}_1$, $\hat{\mP}_2$, $\hat{\mP}_3$, $\hat{\mP}_4$, which are defined as follows:
\begin{equation*}
\begin{aligned}
 \hat{\mP}_1 := \bconv(\{\pbx{1},\pbx{2},\pbx{3},\pbx{4},\pby{4}\}),\ \hat{\mP}_2 := \bconv(\{\pbx{1},\pbx{2},\pbx{3},\pby{3},\pby{4}\}) \\ 
 \hat{\mP}_3 := \bconv(\{\pbx{1},\pbx{2},\pby{2},\pby{3},\pby{4}\}),\ \hat{\mP}_4 := \bconv(\{\pbx{1},\pby{1},\pby{2},\pby{3},\pby{4}\}) 
\end{aligned}
\end{equation*}
A sketch of those can be found in Figure \ref{fig:decomprefprism4}.
\input{figtabs/sec5_prismdecompsketch}
The pentatopes are disjoint (except for a part with measure zero) and sum up to the reference prism. To see this, we give the following characterization of the pentatopes $\hat{\mP}_i$ in terms of constrained sets and their partial sums $\hat{\mB}_i = \bigcup_{j=1}^i \hat{\mP}_j$:
\begin{subequations}
\small
\begin{align*}
  \hat{\mP}_1 &= \{ \bx \in \hat{Q}, x_3 \geq x_4\}, & \!\!\!\! & \\ 
  \hat{\mP}_2 &= \{ \bx \in \hat{Q}, x_3 \leq x_4, x_3+x_2 \geq x_4\}, & \!\!\!\!
  \hat{\mB}_2 &= \{ \bx \in \hat{Q}, x_3+x_2 \geq x_4\}, \\
  \hat{\mP}_3 &= \{ \bx \in \hat{Q}, x_1+x_2+x_3 \geq x_4, x_3+x_2 \leq x_4 \}, & \!\!\!\!
  \hat{\mB}_3 &= \{ \bx \in \hat{Q}, x_1+x_2+x_3 \geq x_4 \}, & \!\!\!\! & \\
  \hat{\mP}_4 &= \{ \bx \in \hat{Q}, x_1+x_2+x_3\leq x_4 \}, & 
  \hat{\mB}_4 &= \{ \bx \in \hat{Q} \}. & \!\!\!\! &
\end{align*}
\end{subequations}

Note further that the measure of all pentatopes are the same, i.e. $\mathrm{meas}_4(\hat{\mP}_i) = 1/24$.

\subsection{Decomposing the reference hypertriangle} \label{sec:decomprefhypertrig}
Let $\bu^i = \hat{\bx}^{1,i}$, $\bv^i = \hat{\bx}^{2,i}$, $\bw^i = \hat{\bx}^{3,i}$, $i=1,..,3$ with $\hat{\bx}^{i,j}$ as in definition \ref{def:hypertriangle}. We decompose $\hat{\mathcal{H}}$ into six pentatopes which are defined as follows:
\begin{equation*}
\begin{array}{l@{}ll}
 \hat{\mD}_u &= \bconv(\{\pbbu1,\pbbu2,\pbbu3,\bv^2,\bw^3\}),\ & \hat{\mD}_v = \bconv(\{\bu^1,\pbbv1,\pbbv2,\pbbv3,\bw^3\}) \\ 
 \hat{\mD}_w &= \bconv(\{\bu^1,\bv^2,\pbbw1,\pbbw2,\pbbw3\}),\ & \hat{\mD}_1 = \bconv(\{\pbbe{u}1,\pbbe{v}1,\bv^2,\pbbe{w}1,\bw^3\}) \\
 \hat{\mD}_2 &= \bconv(\{\bu^1,\pbbf{u}2,\pbbf{v}2,\pbbf{w}2,\bw^3\}),\ & \hat{\mD}_3 = \bconv(\{\bu^1,\pbbg{u}3,\bv^2,\pbbg{v}3,\pbbg{w}3\}) 
\end{array}
\end{equation*}
Note that there is a simple structure behind this decomposition. We define the ``diagonal triangle'' as $\hat{K}_{\mbox{diag}} = \bconv(\bu^1,\bv^2,\bw^3)$. To the three vertices of $\hat{K}_{\mbox{diag}}$ we add the missing vertices (underlined) of one of the following six triangles
\begin{equation*}
\begin{array}{l@{\ \ }l@{\ \ }l}
{\color{mildred} \hat{K}_u }= \bconv(\{\bu^1,\bu^2,\bu^3\}), &
{\color{mildgreen} \hat{K}_v }= \bconv(\{\bv^1,\bv^2,\bv^3\}), &
{\color{mildblue} \hat{K}_w }= \bconv(\{\bw^1,\bw^2,\bw^3\}), \\
{\color{yellow!40!black} \hat{K}_1 }= \bconv(\{\bu^1,\bv^1,\bw^1\}), &
{\color{orange!60!black} \hat{K}_2 }= \bconv(\{\bu^2,\bv^2,\bw^2\}), &
{\color{purple!60!black} \hat{K}_3 }= \bconv(\{\bu^3,\bv^3,\bw^3\}).
\end{array}
\end{equation*}
A sketch of those pentatopes is given in Figure \ref{fig:decomprefhypertrig}.
\input{figtabs/sec5_hypertrigdecompsketch}
Also here, one can show that the pentatopes are disjoint (except for a part with measure zero), and sum up to $\hat{\mH}$. 
To this end we divide the hypertriangle according to three binary decisions and define 
\begin{align*}
\hat{\mH}^{i,j,k} := \hat{\mH} & \cap \ \{ (-1)^i x_2 \leq (-1)^i x_4\} \cap\  \{ (-1)^j x_1 \leq (-1)^j x_3\}  \\ 
 & \cap\  \{ (-1)^k (x_1 + x_2) \leq (-1)^k (x_3 + x_4)\}, \quad i,j,k=0,1 
\end{align*}
Note that $\hat{\mH}^{1,1,0}$ and $\hat{\mH}^{0,0,1}$ are sets with measure zero. All other sets can be identified with a pentatope from the decomposition:
$$
\hat{\mH}^{1,1,1} = \hat{\mD}_u, \ \hat{\mH}^{1,0,0} = \hat{\mD}_v,\  \hat{\mH}^{0,1,0} = \hat{\mD}_w,\ 
\hat{\mH}^{0,0,0} = \hat{\mD}_1, \ \hat{\mH}^{0,1,1} = \hat{\mD}_2,\  \hat{\mH}^{1,0,1} = \hat{\mD}_3. 
$$

\subsection{Decomposition of a pentatope intersected by the space-time interface} \label{sec:decompintersectpenta}
We assume that the space-time interface is approximated in a piecewise planar fashion, s.t. within each pentatope the space-time interface is a (hyper-)plane. This plane divides a pentatope into two parts. Note that due to the pentatope being a convex set each of the two parts will still be convex. 
We now consider a pentatope $\mP$ which is cut by the plane $\mG=\{\bx \in \rr^4: \bx \cdot \bng = c \}$ which represents the local approximation of the space-time interface. Each vertex $\bv$ is marked corresponding to one of the two halfspaces. Vertices with $\bv \cdot \bng < c$ are marked with a plus (+), all others with a minus (-). Note that this classification includes the cases where the space-time interface hits vertices ($\bv \cdot \bng = c$).
We thus can only have two non-trivial situations:
\begin{enumerate}
\item[Case 1:] One vertex has a sign that is different from all the others or
\item[Case 2:] Two vertices have a sign that is different from the other three vertices.
\end{enumerate}
In the following we will consider these cases separately and construct a decomposition of the parts into pentatopes. Without loss of generality we assume that the vertices in the smaller group of vertices are those marked with a plus (+).
\subsubsection{Case 1: Decomposition into one pentatope and one 4-prism} \label{sec:case1}
We consider the case where one vertex of a pentatope, say $\bx^5$, is marked with a plus (+). All other vertices ($\bx^1, \bx^2, \bx^3, \bx^4$) are marked with a minus (-).
The cutting points of the hyperplane $\mG$ with the edges are $\bb^1 := \overline{\bx^1\bx^5} \cap \mG, \bb^2 := \overline{\bx^2\bx^5} \cap \mG, \bb^3 := \overline{\bx^3\bx^5} \cap \mG, \bb^4 := \overline{\bx^4\bx^5} \cap \mG$. 
The geometry containing the separated vertex is the pentatope $\mP^+ := \bconv(\{ \bb^1, \bb^2, \bb^3, \bb^4, \bx^5 \})$ while the remainder is $\mathcal{Q}^- := \bconv(\{ \bx^1, \bx^2, \bx^3, \bx^4 ,\bb^1, \bb^2, \bb^3, \bb^4 \}) $. 
Consider the mapping
\vspace*{-0.1cm}
$$
\Phi: \hat{\mQ} \rightarrow \mQ^-, (\hat{x}_1,\hat{x}_2,\hat{x}_3,\hat{x}_4) \rightarrow \sum_{i=1}^4 \mu_i(\hat{x}_1,\hat{x}_2,\hat{x}_3)(\hat{x}_4 \bb^i+(1-\hat{x}_4)\bx^i) 
\vspace*{-0.2cm}
$$
with $\hat{\mu}(\hat{x}_1,\hat{x}_2,\hat{x}_3)$ the barycentric coordinates of the reference tetrahedron $\hat{T}$.
The decomposition of the reference 4-prism $\hat{\mQ}$ into the four pentatopes $\hat{\mP}_i,\ i=1,..,4$ as described in section \ref{sec:decomprefprism} can be used as a triangulation of $\hat{\mQ}$. Let $\Phi_h$ be the (pentatope-) piecewise linear interpolation of $\Phi$ at the vertices of this triangulation. 
One can show that $\Phi_h$ is an isomorphism between $\hat{\mQ}$ and $\mQ^-$ and furthermore for each pentatope $\hat{\mP}_i$ the image $\Phi_h(\hat{\mP}_i)$ is again a pentatope. 
Thus the decomposition rule for the reference 4-prism can also be applied here and we get a valid decomposition by taking the four pentatopes
\begin{equation*}
\begin{aligned}
 \mP_1 \! = \! \Phi_h(\hat{\mP}_1) =  \bconv(\{\pbx{1},\pbx{2},\pbx{3},\pbx{4},\pbb{4}\}),\ \mP_2 \! = \! \Phi_h(\hat{\mP}_2) = \bconv(\{\pbx{1},\pbx{2},\pbx{3},\pbb{3},\pbb{4}\}) \\ 
 \mP_3 \! = \! \Phi_h(\hat{\mP}_3) = \bconv(\{\pbx{1},\pbx{2},\pbb{2},\pbb{3},\pbb{4}\}),\ \mP_4 \! = \! \Phi_h(\hat{\mP}_4) = \bconv(\{\pbx{1},\pbb{1},\pbb{2},\pbb{3},\pbb{4}\}). 
\end{aligned}
\end{equation*}

\paragraph{Decomposition of the space-time interface into tetrahedra for case 1}
The triangulation of the interface is trivially obtained with the tetrahedron
$$ 
\mP \cap \mG = \mI = \bconv (\{ \bb^1, \bb^2, \bb^3, \bb^4\}).
$$
\subsubsection{Case 2: Decomposition into one 4-prism and one hypertriangle} \label{sec:case2}
Let us consider the case where two vertices of a pentatope are marked with a plus (+), these are (w.l.o.g.) vertices $\bx^4$ and $\bx^5$. All other vertices ($\bx^1, \bx^2, \bx^3$) are marked with a minus (-).
The cutting points of the hyperplane $\mG$ with the edges are 
$\bc^1 := \overline{\bx^1\bx^4} \cap \mG$,
$\bc^2 := \overline{\bx^2\bx^4} \cap \mG$,
$\bc^3 := \overline{\bx^3\bx^4} \cap \mG$,
$\bd^1 := \overline{\bx^1\bx^5} \cap \mG$,
$\bd^2 := \overline{\bx^2\bx^5} \cap \mG$,
$\bd^3 := \overline{\bx^3\bx^5} \cap \mG$.
Thus we have to decompose the two parts $\mathcal{H}^-$ and $\mathcal{Q}^+$ into pentatopes with
\begin{eqnarray*}
\mathcal{H}^-&:=&\bconv(\{\bx^1,\bx^2,\bx^3,\bc^1,\bc^2,\bc^3,\bd^1,\bd^2,\bd^3\}),\\
\mathcal{Q}^+&:=&\bconv(\{\bc^1,\bc^2,\bc^3,\bd^1,\bd^2,\bd^3,\bx^4,\bx^5\}).
\end{eqnarray*}
Let us start with the decomposition of $\mathcal{H}^-$. 
Consider the mapping
\vspace*{-0.1cm}
$$ 
\Phi: \hat{\mH} \rightarrow \mH^-, (\hat{x}_1,\hat{x}_2,\hat{x}_3,\hat{x}_4) \rightarrow \sum_{i=1}^3 \sum_{j=1}^3 \rho_i(\hat{x}_1,\hat{x}_2) \rho_j(\hat{x}_3,\hat{x}_4) \bq^{i,j} 
\vspace*{-0.2cm}
$$
with $\bq^{i,1} = \bx^i$, $\bq^{i,2} = \bc^i$ and $\bq^{i,3} = \bd^i$
where $\rho_i(\hat{x}_1,\hat{x}_2) $ are the barycentric coordinates of the reference triangle $\hat{K} \subset \rr^2$. Following section \ref{sec:decomprefhypertrig}, we have a triangulation of $\hat{\mH}$ into pentatopes  $\{\hat{\mD}_i\}$. One can show that the (pentatope-) piecewise linear interpolation $\Phi_h$ of $\Phi$ is an isomorphism between $\hat{\mH}$ and $\mH^-$ and each image $\Phi_h(\hat{\mD}_i)$ is again a pentatope. Therefore we can apply the decomposition of the reference hypertriangle $\hat{\mH}$ into pentatopes to get the six pentatopes
\begin{eqnarray*}
\mD_u=\Phi_h(\hat{\mD}_u) &=& \bconv(\{\pbbx1,\pbbx2,\pbbx3,\bc^2,\bd^3\}), \\
\mD_v=\Phi_h(\hat{\mD}_v) &=& \bconv(\{\bx^1,\pbbc1,\pbbc2,\pbbc3,\bd^3\}), \\ 
\mD_w=\Phi_h(\hat{\mD}_w) &=& \bconv(\{\bx^1,\bc^2,\pbbd1,\pbbd2,\pbbd3\}), \\
\mD_1=\Phi_h(\hat{\mD}_1) &=& \bconv(\{\pbbe{x}1,\pbbe{c}1,\bc^2,\pbbe{d}1,\bd^3\}),\\
\mD_2=\Phi_h(\hat{\mD}_2) &=& \bconv(\{\bx^1,\pbbf{x}2,\pbbf{c}2,\pbbf{d}2,\bd^3\}),\\ 
\mD_3=\Phi_h(\hat{\mD}_3) &=& \bconv(\{\bx^1,\pbbg{x}3,\bc^2,\pbbg{c}3,\pbbg{d}3\}.
\end{eqnarray*}
We now turn over to $\mathcal{Q}^+$. For notational convenience define $\pbc{4} := \pbx{4}$ and $\pbd{4} := \pbxx{5}$. Thus  $\mathcal{Q}^+ = \bconv(\{\pbc{1},\pbc{2},\pbc{3},\pbc{4},\pbd{1},\pbd{2},\pbd{3},\pbd{4}\})$. Now the structure is similar to the situation for $\mQ^-$ in Case 1 and we can apply the same procedure and get a valid decomposition $\bigcup \mP_i = \mathcal{Q}^+$ with
\begin{equation*}
\begin{aligned}
 \mP_1 \! = \! \Phi_h^\mQ(\hat{\mP}_1) =  \bconv(\{\pbc{1},\pbc{2},\pbc{3},\pbx{4},\pbxx{5}\}),\ \mP_2 \! = \! \Phi_h^\mQ(\hat{\mP}_2) = \bconv(\{\pbc{1},\pbc{2},\pbc{3},\pbd{3},\pbxx{5}\}) \\ 
 \mP_3 \! = \! \Phi_h^\mQ(\hat{\mP}_3) = \bconv(\{\pbc{1},\pbc{2},\pbd{2},\pbd{3},\pbxx{5}\}),\ \mP_4 \! = \! \Phi_h^\mQ(\hat{\mP}_4) = \bconv(\{\pbc{1},\pbd{1},\pbd{2},\pbd{3},\pbxx{5}\})
\end{aligned}
\end{equation*}
with $\Phi_h^\mQ$ the corresponding piecewise linear transformation for the 4-prism.
\paragraph{Decomposition of the space-time interface into tetrahedra for case 2}
With similar techniques as done for the four dimensional volume, we can proceed with the triangulation of the interface which is isomorph to a 3-prism resulting in tetrahedra $\mI_i$:
\begin{equation*}
\mI_1 \! =\!  \bconv(\{ \pbc{1}, \pbc{2}, \pbc{3}, \pbd{3}\}), \ 
\mI_2 \! =\!  \bconv(\{ \pbc{1}, \pbc{2}, \pbd{2}, \pbd{3}\}), \ 
\mI_3 \! =\!  \bconv(\{ \pbc{1}, \pbd{1}, \pbd{2}, \pbd{3}\}).
\end{equation*} 


\section*{Conclusion}
The method presented and analysed in \cite{ReuskenLehrenfeld13} is applied to a spatially three dimensional mass transport problem with a moving discontinuity. A strategy for constructing a piecewise planar approximate space-time interface and corresponding polygonal space-time domains is presented. For the integration on the polygonal space-time domains a method which decomposes the domain into a set of simplices is presented. Furthermore the construction of the tensor product space-time finite elements and the (non-standard) issue of quadrature rules in four space dimensions is addressed. Results of numerical experiments reveal at least a second order convergence in time.

\subsection*{Acknowledgement}
The author wants to thank Arnold Reusken and the referees for their comments that greatly improved the manuscript.
Financial support from the German Science Foundation (DFG) within the Priority Program (SPP) 1506 ``Transport Processes at Fluidic Interfaces'' is gratefully acknowledged.


\appendix
\section{Quadrature rule(s) on pentatopes}\label{sec:quad4d}
We shortly review lower order quadrature rules on pentatopes and discuss how to achieve higher order rules using tetrahedron rules, Duffy transformation and 1D Gauss-Jacobi integration rules. Integration rules are given for the reference pentatope $\hat{\mP}$.
\subsection{First order rule}
There holds  $\int_{\hat{\mP}} 1 \, d\hat{\bx} = 1/24$ and  $\int_{\hat{\mP}} q(\hat{\bx}) \, d\hat{\bx} = 1/120$ 
for $q(\hat{\bx})\in \{\hat{x}_1,\hat{x}_2,\hat{x}_3,\hat{x}_4\}$, s.t. 
the following rule is obviously exact for all polynomials up to degree one:
\vspace*{-0.35cm} 
$$ 
I^1(f) = 1/120 \sum_{i=0}^4 f(\be^i) 
\vspace*{-0.1cm} 
$$ 
\subsection{Third order rule}\label{sec:quad4d:thirdorder}
A third order rule, taken from \cite{Behr08}, is as follows:
\vspace*{-0.1cm} 
$$ I^3(f) = 1/120 \sum_{i=1}^5 f(\bs^i) \mbox{ with } \lambda_j(\bs^i) = \alpha \mbox{ for } i\neq j \mbox{ and } \lambda_j(\bs^i) = \beta \mbox{ for } i=j 
\vspace*{-0.1cm} 
 $$ 
where $\lambda_j$ is the barycentric coordinate of vertex $j$ in the reference pentatope and the coefficients are $\displaystyle \alpha = 0.118350341907227374 \mbox{ and } \beta = 0.526598632371090503$.
\subsection{Higher order rules using the Duffy transformation}
\begin{figure}
  \centering
  \begin{tikzpicture}
    [
    rotate=0,
    ]
    \begin{scope}[scale=1.1]
      \draw[->,very thick,>=stealth'] (-0.1,0)  -- (2.5,0) node(xline)[below]{$y$};
      \draw[->,very thick,>=stealth'] (0,-0.1)  -- (0,2.5) node(yline)[left]{$t$};
      \draw [] (2.0,0) -- (0,2.0); 
      \draw [] (2.0,-0.1) -- (2.0,0.1);
      \draw [] (2.0,-0.1) node[below]{1};
      \draw [] (-0.1,2.0) -- (0.1,2.0);
      \draw [] (-0.1,2.0) node[left]{1};
      \node[anchor=center] at (1.0, -0.65) {$0 \leq t \leq 1, \ \ 0 \leq y \leq 1 -t$};
    \end{scope}
    \draw [<->] (3.1,1.5) -- (5.9,1.5);
    \node[above]  at (4.5,1.5) {$\begin{array}{c@{}c@{}c} \tilde{t} &= & t,\\ \tilde{y} &=& \frac{1}{1-t} y\end{array}$};
    \node[below]  at (4.5,1.5) {$\begin{array}{c@{}c@{}c} t &= & \tilde{t},\\ y &=& (1-\tilde{t}) \tilde{y}\end{array}$};
    \begin{scope}[xshift=6.3cm,scale=1.1]
      \draw[->,very thick,>=stealth'] (-0.1,0)  -- (2.5,0) node(xline)[below]{$\tilde{y}$};
      \draw[->,very thick,>=stealth'] (0,-0.1)  -- (0,2.5) node(yline)[left]{$\tilde{t}$};
      \draw [] (2.0,0) -- (2.0,2.0); 
      \draw [] (2.0,2.0) -- (0,2.0); 
      \draw [] (2.0,-0.1) -- (2.0,0.1);
      \draw [] (2.0,-0.1) node[below]{1};
      \draw [] (-0.1,2.0) -- (0.1,2.0);
      \draw [] (-0.1,2.0) node[left]{1};
      \node[anchor=center] at (1.0, -0.65) {$0 \leq \tilde{t} \leq 1, \ \ 0 \leq \tilde{y} \leq 1$};
    \end{scope}
  \end{tikzpicture}
  \vspace*{-0.3cm}
  \caption{Sketch of the Duffy transformation for $d=1$}
  \label{fig:duffy}
\vspace*{-0.25cm}
\end{figure}
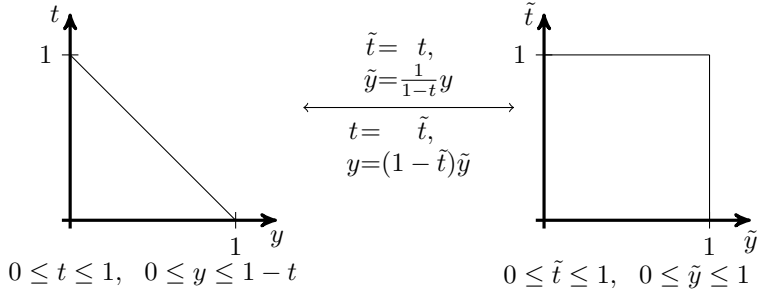

A more general approach to derive integration rules for pentatopes is based on the Duffy transformation \cite{duffy82}.
Let $\hat{\by} = (\hat{y}_1,\hat{y}_2,\hat{y}_3) \in \rr^3$, and $\hat{\bx} = (\hat{\by},t) \in \rr^4$.
The problem to compute $ \int_{\hat{\mP}} f(\hat{\bx}) \, d \hat{\bx} = \int_{\hat{\mP}} f(\hat{\by},t) \, d (\hat{\by},t) $ can be transformed using the transformation $ (\hat{\by},t) \rightarrow ( 1/(1-t) \hat{\by},t) = (\tilde{\by},t) $ (see also Figure \ref{fig:duffy} for a sketch):
\begin{eqnarray*}
  && \int_{\hat{\mP}} f(\hat{\by},t) \, d (\hat{\by},t) =  \int_0^1 
  \int_0^{1-t}
  \int_0^{1-t-\hat{y}_1}
  \int_0^{1-t-\hat{y}_1-\hat{y}_2}
  \!\!\!\!\!\!\!\!\!\!\!\!\!\!
  f(\hat{\by},t) \, d\hat{y}_3 \, d\hat{y}_2 \, d\hat{y}_1  \, dt \\
  && {\text{\small [$\tilde{\by} = 1/(1-t) \hat{\by}$]}\qquad }  \\
  && =  \int_0^1 
  (1-t)^3
  \int_0^1
  \int_0^{1-\tilde{y}_1}
  \!\!\!\!
  \int_0^{1-\tilde{y}_1-\tilde{y}_2}
  \!\!\!\!\!\!\!\!\!\! f( (1-t)\tilde{\by},t) \, d\tilde{y}_3 \, d\tilde{y}_2 \, d\tilde{y}_1 \, dt \\
  && = \int_0^1 (1-t)^3 \int_{\hat{T}} \tilde{f}(\tilde{\by},t) \, d \tilde{\by} \, dt
  =
  \int_0^1 (1-t)^3 \tilde{g}(t)\, dt 
\end{eqnarray*}
with $\tilde{f}(\tilde{\by},t) = f((1-t)\tilde{\by},t)$ and $\tilde{g}(t) = \int_{\hat{T}} \tilde{f}(\tilde{\by},t) \, d \tilde{\by}$ . 
In this form one can apply a one-dimensional integration rule of the form
\begin{eqnarray*}
  \int_0^1 (1-t)^3 \tilde{g}(t)\, dt \approx \sum_{k=0}^{N} \omega_i \tilde{g}(t_i)
\end{eqnarray*}
where $\omega_i$ and $t_i$ are weights and points of the corresponding quadrature rule. In order to approximate $\tilde{g}(t_i)$ at every integration point $t_i$ a standard 3D quadrature rule can be applied. Let's assume this 3D quadrature rule has order $q$ accuracy. The highest order for the pentatope rule at lowest costs is achieved if a Gauss-Jacobi rule (corresponding to the weight $(1-t)^3$) of order $q$ is used for the numerical integration w.r.t. $t$.  The resulting quadrature rule is positive, but not symmetric. In principle also the quadrature rule for the tetrahedron can be derived from lower dimensional quadrature rules applying the idea recursively. 
This generic procedure generates quadrature rules which have slightly more points than symmetric Gauss rules.  
For that reason, we use symmetric Gauss rules for the tetrahedron. 
\section{Computing the weighting factor $\nu$}\label{sec:nu}
The weighting factor $\nu(\bs)$ in the Nitsche XFEM-DG method can be computed using the space-time normal $\bn_{\Gamma_\ast}$ of the space-time interface. One can show that there holds 
$$\nu(\bs) = (1+(\bw \cdot \bn_\Gamma)^2)^{-\frac12} = \Vert (n_1,..,n_{d})^T \Vert, \quad \bs \in \Gs $$
with $\bn_{\Gamma_\ast} = (n_1,..,n_{d+1})^T$ the space-time normal at the interface.

As we use a piecewise planar approximation of the space-time interface consisting of d-simplices in $d+1$ dimensions we have to compute a normal to the d-simplex. It is known that for $d=2$ one can use the standard cross-product to compute the normal. In the next section we quote a generalized cross-product which allows to do the same if $d=3$.
\subsection{Computing normals to tetrahedra in 4 dimensions}\label{sec:normals4D}
In \cite{Holasch} a generalization of the cross-product is given. Given three vectors $\bu^1, \bu^2, \bu^3 \in \rr^4$ one can compute the cross-product $\bv = \bX(\bu^1, \bu^2, \bu^3)$, s.t.
\begin{itemize}
\item $\bX(\bu^1, \bu^2, \bu^3)=0$ iff $\bu^1, \bu^2, \bu^3$ are linear dependent.
\item Iff $\bu^1, \bu^2, \bu^3$ are linear independent then for $\bv = \bX(\bu^1, \bu^2, \bu^3)$, there holds: $\bv \perp \bu^i$, $i=1,..,3$.
\item $\alpha \bX(\bu^1, \bu^2, \bu^3)= \bX( \alpha \bu^1, \bu^2, \bu^3) = \bX( \bu^1, \alpha \bu^2, \bu^3) = \bX( \bu^1, \bu^2, \alpha \bu^3)$, $\alpha \in \rr$
\item $\bX(\bu^1, \bu^2, \bu^3) = \mathrm{sign}(\pi) \bX(\bu^{\pi(1)}, \bu^{\pi(2)}, \bu^{\pi(3)})$, where $\pi$ is a permutation, i.e. changing the order of the arguments switches the sign.  
\end{itemize}
This cross-product can be used to compute normals to tetrahedra.  The computation is given below:

Given $\bu, \bv, \bw \in \rr^4$. Compute $\bz = \bX(\bu, \bv, \bw) \in \rr^4$ as follows:
{
\vspace*{-0.2cm}
\begin{equation*}
  \begin{array}{c@{}c}
    a_{1,2} &= u_1 \cdot v_2 - u_2 \cdot v_1, \\
    a_{1,3} &= u_1 \cdot v_3 - u_3 \cdot v_1, \\
    a_{1,4} &= u_1 \cdot v_4 - u_4 \cdot v_1, \\
    a_{2,3} &= u_2 \cdot v_3 - u_3 \cdot v_2, \\
    a_{2,4} &= u_2 \cdot v_4 - u_4 \cdot v_2, \\
    a_{3,4} &= u_3 \cdot v_4 - u_4 \cdot v_3, \\
  \end{array}
  \quad\quad
  \begin{array}{c@{}l@{}c}
    z_1 &=   & w_2 \cdot a_{3,4} - w_3 \cdot a_{2,4} + w_4 \cdot a_{2,3} \\
    z_2 &= - & w_1 \cdot a_{3,4} + w_3 \cdot a_{1,4} - w_4 \cdot a_{1,3} \\
    z_3 &=   & w_1 \cdot a_{2,4} - w_2 \cdot a_{1,4} + w_4 \cdot a_{1,2} \\
    z_4 &= - & w_1 \cdot a_{2,3} + w_2 \cdot a_{1,3} - w_3 \cdot a_{1,2} \\
  \end{array}
\vspace*{-0.3cm}
\end{equation*}
}
\vspace*{-0.3cm}

\bibliographystyle{siam}
\bibliography{literature}

\end{document}